\documentclass{amsart}
\usepackage{amsmath, amscd, amsthm, amssymb, amsfonts, latexsym, stmaryrd}

\usepackage{geometry}
\usepackage{graphicx, float, color}
\usepackage[all]{xy}
\usepackage{caption, subcaption}

\usepackage{tikz}
\usepackage{pgfplots}
\pgfplotsset{compat=1.18}
\usetikzlibrary{positioning, quotes, calc, arrows.meta, shapes.geometric}

\tikzset{
    base/.style={line width=0.6pt, >=Stealth},
    dot/.style={circle, fill, inner sep=1.2pt},
    every edge quotes/.style={auto, font=\small}
}

\usepackage{url}
\usepackage{hyperref}
\hypersetup{
    colorlinks=false,
    linkcolor=blue,
    filecolor=magenta,      
    urlcolor=cyan,
    pdftitle={Character varieties on a four-holed sphere},
    pdfpagemode=FullScreen,
}

\numberwithin{equation}{section}

\def\Bigroman{\uppercase\expandafter{\romannumeral\number\count 255 }}
\def\Romannumeral{\afterassignment\BigRoman\count255=}

\theoremstyle{definition}
\newtheorem{definition}{Definition}[section]

\theoremstyle{remark}
\newtheorem{remark}{Remark}[section]
\newtheorem{example}[definition]{Example}

\theoremstyle{plain}
\newtheorem{theorem}{Theorem}[section] 

\newtheorem{lemma}[theorem]{Lemma}

\newtheorem{claim}[theorem]{Claim}

\DeclareMathOperator{\SL}{SL}
\DeclareMathOperator{\tr}{tr}
\DeclareMathOperator{\Hom}{Hom}
\DeclareMathOperator{\Diff}{Diff}
\DeclareMathOperator{\MCG}{MCG}

\title{Character varieties on a four-holed sphere}
\author{Eunju Shin}
\address{Department of Mathematics, Seoul National University, 1 Gwanak-ro, Gwanak-gu, Seoul 08826\\ Current address: Department of Mathematics, The Ohio State University, 231 West 18th Avenue, Columbus, OH 43210-1174}
\email{shin.970@osu.edu}
\thanks{This is the author's accepted manuscript of an article accepted for publication in \textit{Proceedings of the American Mathematical Society}. The final published version will be available at https://doi.org/10.1090/proc/17770}
\usepackage{lipsum}

\begin{document}
\begin{abstract}
    For each $\mathbf{k}\in\mathbb{C}^4$, let $V_\mathbf{k}$ be the character variety on a four-holed sphere and $\Gamma_\mathbf{k}$ the group generated by the Vieta involution maps. First, under a certain condition, we find a fundamental domain for $\Gamma_\mathbf{k}$-action on $V_\mathbf{k}$, expressed via inequalities. Second, we show that it is decidable whether or not two integral solutions for $V_\mathbf{k}$ are in the same $\Gamma_\mathbf{k}$-orbit and in the same mapping class group orbit. To achieve our goals, we introduce graphs corresponding to the $\Gamma_\mathbf{k}$-orbits and the mapping class group orbits, and classify their restricted global shapes by analyzing the limited local edge configurations at each vertex, using a descent argument.
\end{abstract}
\maketitle

\section{Introduction} 
For each $\mathbf{k}=(k_1,k_2,k_3,k_4)\in \mathbb{C}^4$, we write
\begin{align*}
   \begin{cases}
       \alpha_1=k_1k_2+k_3k_4,\\
       \alpha_2=k_1k_4+k_2k_3,\\
       \alpha_3=k_1k_3+k_2k_4,\\
       \beta=4-\sum_{i=1}^4k_i^2-\prod_{i=1}^4k_i.
   \end{cases}
   \end{align*}
   Let $V_\mathbf{k}(\mathbb{Z})$ denote the set of integral solutions for a generalized Markoff equation, given by
   \begin{align*}
       V_\mathbf{k}:x_1^2+x_2^2+x_3^2+x_1x_2x_3-\alpha_1x_1-\alpha_2x_2-\alpha_3x_3-\beta=0.
   \end{align*}
   Let $\Gamma_\mathbf{k}$ be the group generated by three Vieta involution maps
   \begin{align*}
    \mathcal{V}_1(x_1,x_2,x_3)&=(\alpha_1-x_2x_3-x_1,x_2,x_3),\\
    \mathcal{V}_2(x_1,x_2,x_3)&=(x_1,\alpha_2-x_3x_1-x_2,x_3),\\
    \mathcal{V}_3(x_1,x_2,x_3)&=(x_1,x_2,\alpha_3-x_1x_2-x_3).
\end{align*}
   From the equation
\begin{align*}
    &x_1^2+x_2^2+x_3^2+x_1x_2x_3-\alpha_1x_1-\alpha_2x_2-\alpha_3x_3-\beta\\
    &=-x_1(\alpha_1-x_2x_3-x_1)+x_2^2+x_3^2-\alpha_2x_2-\alpha_3x_3-\beta,
\end{align*}
we have that $\mathcal{V}_i$'s send $V_\mathbf{k}(\mathbb{Z})$ to itself. We define an action of $\Gamma_\mathbf{k}$ on $V_\mathbf{k}(\mathbb{Z})$ to be $\gamma\cdot\mathbf{x}:=\gamma(\mathbf{x})$ for $\gamma\in\Gamma_\mathbf{k}$ and $\mathbf{x}\in V_\mathbf{k}(\mathbb{Z})$. (We say $\mathbf{x},\mathbf{y}\in V_\mathbf{k}$ are \emph{$\Gamma_\mathbf{k}$-equivalent} if there exists $\gamma\in\Gamma_\mathbf{k}$ such that $\mathbf{x}=\gamma\cdot\mathbf{y}$.) Moreover, the mapping class group of a four-holed sphere acts on $V_\mathbf{k}(\mathbb{Z})$ via polynomial transformations, given by compositions of pairs of Vieta involution maps $\mathcal{V}_i$;
\begin{align*}
    \mathcal{V}_3\mathcal{V}_2(x_1,x_2,x_3)&=(x_1,\alpha_2-x_3x_1-x_2,\alpha_3-x_1(\alpha_2-x_3x_1-x_2)-x_3),\\
    \mathcal{V}_1\mathcal{V}_3(x_1,x_2,x_3)&=(\alpha_1-x_2(\alpha_3-x_1x_2-x_3)-x_1,x_2,\alpha_3-x_1x_2-x_3),\\
    \mathcal{V}_2\mathcal{V}_1(x_1,x_2,x_3)&=(\alpha_1-x_2x_3-x_1,\alpha_2-x_3(\alpha_1-x_2x_3-x_1)-x_2,x_3).
\end{align*}

The generalized Markoff equation can be interpreted as the character variety on a four-holed sphere (as we will see in Subsection 1.3). This paper deals with finding a fundamental domain for a group action of the group $\Gamma_\mathbf{k}$ generated by the Vieta involution maps $\mathcal{V}_i$ on the set $V_k(\mathbb{Z})$ of integral solutions for the generalized Markoff equation. In general, even for a linear reduction, it is challenging to express a fundamental domain in the set of integral solutions for an equation with at least $3$ variables as a reduced form that satisfies some inequalities\cite{so}. This paper examines the decidability (which will be defined below) of the equivalence of solutions in $V_\mathbf{k}(\mathbb{Z})$ under the Vieta involution action (via compositions of $\mathcal{V}_i$'s) and the mapping class group action.

We start with the motivation. Binary quadratic forms historically have played an essential role in number theory. One classical problem focuses on finding representatives of the equivalence classes of binary quadratic forms with the smallest possible coefficients. It is called reduction theory;
one of the most famous examples in this theory is Lagrange-Gauss reduction. Grunewald and Segal\cite{gs2} constructed an algorithm to decide whether a given homogeneous space for an arithmetic group has an integral point, as a reduction for arithmetic groups. In this paper, we deal with a nonlinear analogue of the Grunewald-Segal algorithm; to be precise, we establish a nonlinear reduction on the generalized Markoff type cubic affine surfaces for $\Gamma_\mathbf{k}$-action. It is worth mentioning that in the specific case of $\alpha_1=\alpha_2=\alpha_3=0$ and $\beta=k\neq4$ for the generalized Markoff equation, Ghosh and Sarnak\cite{gs} recently found a fundamental domain for $\Gamma_\mathbf{k}$-action on the set of integral solutions for $x^2+y^2+z^2+xyz=k$, expressed in terms of inequalities.

\subsection{Main results}
We describe now the main results of this paper.

\begin{definition}[fundamental domain] Let $\mathcal{S}$ be a set, and let a group $\Gamma$ act on $\mathcal{S}$. If a subset $\mathcal{F}\subset\mathcal{S}$ satisfies the following conditions:
\begin{itemize}
    \item[(i)] $\mathcal{S}=\bigcup_{\gamma\in\Gamma}\gamma\cdot\mathcal{F}$, where $\gamma\cdot\mathcal{F}=\{\gamma\cdot f:f\in\mathcal{F}\}$;
    \item[(ii)] For $\gamma_1, \gamma_2\in\Gamma$, either $\gamma_1\cdot\mathcal{F}=\gamma_2\cdot\mathcal{F}$ or $(\gamma_1\cdot\mathcal{F})\cap(\gamma_2\cdot\mathcal{F})=\emptyset$,
\end{itemize}
then $\mathcal{F}$ is a \emph{fundamental domain} for $\Gamma$ on $\mathcal{S}$.
\end{definition}

Let $\mathcal{I}=\{(1,2,3),(1,3,2),(2,1,3),(2,3,1),(3,1,2),(3,2,1)\}$. We will use this notation throughout the paper. 

To construct the sets appearing in Theorem \ref{1.1}, we collect minimal vertices $\mathbf{x}=(x_1,x_2,x_3)$ satisfying $\min\{|x_1|,|x_2|,|x_3|\}=s$ for each $s=0,1$, or $s=((-\infty,-2]\cup[2,\infty))\cap\mathbb{N}$. Each such set is then subdivided to ensure that $\mathbf{x}$ does not satisfy $-x_i=\alpha_i-x_jx_k-x_i$ for any $(i,j,k)\in\mathcal{I}$, which is achieved by choosing, between $\mathbf{x}$ and $\mathcal{V}_i(\mathbf{x})$, the one whose $i$-th coordinate is nonnegative.
(Here, minimal vertices will be defined later.)

The main results of this paper are as follows: 
\begin{theorem}\label{1.1}
    For each $\mathbf{k}\in\mathbb{C}^4$, if $\mathbf{x}=(x_1,x_2,x_3)\in V_\mathbf{k}(\mathbb{Z})$ then $\mathbf{x}$ is $\Gamma_\mathbf{k}$-equivalent to an element in $\mathfrak{S}(\mathbf{k})\cup\mathfrak{T}_0(\mathbf{k})\cup\mathfrak{T}_{1}(\mathbf{k})\cup\mathfrak{T}_{-1}(\mathbf{k})$, where
    \begin{align*}
        \mathfrak{S}(\mathbf{k})&=\left\{\mathbf{x}\in V_\mathbf{k}(\mathbb{Z}):\begin{matrix}
            \forall (i,j,k)\in\mathcal{I}\text{ with } x_jx_k\neq\alpha_i,\ 2\leq|x_i|<|\alpha_i-x_jx_k-x_i|\\
            \text{and }\forall (s,t,r)\in\mathcal{I}\text{ with } x_tx_r=\alpha_s,\ 2\leq|x_s|=x_s
        \end{matrix}\right\},\\
        \mathfrak{T}_0(\mathbf{k})&=\left\{\mathbf{x}\in V_\mathbf{k}(\mathbb{Z}):
                \exists s\in\{1,2,3\}\text{ s.t. }x_s=0\
            \&\ \forall t\neq s,\ |x_t|=\min\{|x_t|,|\alpha_t-x_t|\}\ \&\ \alpha_t\neq0\right\}\\
            &\hspace{13pt}\cup\left\{\mathbf{x}\in V_\mathbf{k}(\mathbb{Z}):
            \begin{matrix}
                \exists (s,t,r)\in\mathcal{I}\text{ s.t. }x_s=0\ \&\ |x_t|=x_t\\
            \&\ |x_r|=\min\{|x_r|,|\alpha_r-x_r|\}\ \&\ \alpha_t=0\ \&\ \alpha_r\neq0
            \end{matrix}\right\}\\
            &\hspace{13pt}\cup\left\{\mathbf{x}\in V_\mathbf{k}(\mathbb{Z}):
            \exists s\in\{1,2,3\}\text{ s.t. } x_s=0\ \&\ \forall t\neq s,\ |x_t|=x_t\ \&\ \alpha_t=0\right\},\\
        \mathfrak{T}_{1}(\mathbf{k})&=\left\{\mathbf{x}\in V_\mathbf{k}(\mathbb{Z}):\begin{matrix}\exists s\in\{1,2,3\}\text{ s.t. }\forall(s,t,r)\in\mathcal{I},\ x_t\neq\alpha_r\ \\
            \&\ x_s=1\leq\min\left\{\begin{matrix}
                |\alpha_s-x_tx_r-1|,|\alpha_s-(\alpha_t-x_r-x_t)x_r-1|,\\
                |\alpha_s-(\alpha_t-x_r-x_t)(\alpha_r-\alpha_t+x_t)-1|,\\
                |\alpha_s-(\alpha_t-\alpha_r+x_r)(\alpha_r-\alpha_t+x_t)-1|
            \end{matrix}\right\}\\
            \&\ 1\leq|x_t|=\min\{|x_t|,|\alpha_t-x_r-x_t|,|\alpha_t-\alpha_r-x_r|\}
        \end{matrix}\right\}\\
        &\hspace{13pt}\cup\left\{\mathbf{x}\in V_\mathbf{k}(\mathbb{Z}):\begin{matrix}\exists (s,t,r)\in\mathcal{I}\text{ s.t. }x_t\neq\alpha_r\ \&\  x_r=\alpha_t\\
            \&\ x_s=1\leq\min\left\{\begin{matrix}
                |\alpha_s\pm\alpha_tx_t-1|,|\alpha_s- x_t(x_t\pm(\alpha_r-\alpha_t))-1|,\\
                |\alpha_s\pm (2\alpha_t-\alpha_r)(\alpha_r-\alpha_t\pm x_t)-1|
            \end{matrix}\right\}\\
            \&\ 1\leq|\alpha_t|=\min\{|\alpha_t|,|\alpha_r-\alpha_t\pm x_t|\}\
            \&\ 1\leq x_t=\min\{|x_t|,|2\alpha_t-\alpha_r|\}
        \end{matrix}\right\}\\
        &\hspace{13pt}\cup\left\{\mathbf{x}\in V_\mathbf{k}(\mathbb{Z}):
            \exists s\in\{1,2,3\}\text{ s.t. } x_s=1\ \&\ \forall (s,t,r)\in\mathcal{I},\ |x_t|=x_t=\alpha_r\right\},\\
        \mathfrak{T}_{-1}(\mathbf{k})&=\left\{\mathbf{x}\in V_\mathbf{k}(\mathbb{Z}):\begin{matrix}\forall i\in\{1,2,3\}\text{ with }x_i=-1,\ \alpha_i\neq x_jx_k\text{ for }(i,j,k)\in\mathcal{I}\\
        \text{and }\exists s\in\{1,2,3\}\text{ s.t. }\forall(s,t,r)\in\mathcal{I},\ x_t\neq-\alpha_r\\
            \&\ -x_s=1\leq\min\left\{\begin{matrix}
                |\alpha_s-x_tx_r+1|,|\alpha_s-(\alpha_t+x_r-x_t)x_r+1|,\\
                |\alpha_s-(\alpha_t+x_r-x_t)(\alpha_r+\alpha_t-x_t)+1|,\\
                |\alpha_s-(\alpha_t+\alpha_r-x_r)(\alpha_r+\alpha_t-x_t)+1|
            \end{matrix}\right\}\\
            \&\ 1\leq|x_t|=\min\{|x_t|,|\alpha_t+x_r-x_t|,|\alpha_t+\alpha_r-x_r|\}
        \end{matrix}\right\}\\
        &\hspace{13pt}\cup\left\{\mathbf{x}\in V_\mathbf{k}(\mathbb{Z}):\begin{matrix}\forall i\in\{1,2,3\}\text{ with }x_i=-1,\ \alpha_i\neq x_jx_k\text{ for }(i,j,k)\in\mathcal{I}\\
        \text{and }\exists (s,t,r)\in\mathcal{I}\text{ s.t. } x_t=-\alpha_r\ \&\ x_r\neq-\alpha_t\\
            \&\ -x_s=1\leq\min\left\{\begin{matrix}
                |\alpha_s\pm\alpha_tx_t+1|,|\alpha_s- x_t(x_t\pm(\alpha_r+\alpha_t))+1|,\\
                |\alpha_s\pm (2\alpha_t+\alpha_r)(\alpha_r+\alpha_t\pm x_t)+1|
            \end{matrix}\right\}\\
            \&\ 1\leq|\alpha_t|=\min\{|\alpha_t|,|\alpha_r+\alpha_t\pm x_t|\}\
            \&\ 1\leq x_t=\min\{|x_t|,|2\alpha_t+\alpha_r|\}
        \end{matrix}\right\}\\
        &\hspace{13pt}\cup\left\{\mathbf{x}\in V_\mathbf{k}(\mathbb{Z}):\begin{matrix}
            \forall i\in\{1,2,3\}\text{ with }x_i=-1,\ \alpha_i\neq x_jx_k\text{ for }(i,j,k)\in\mathcal{I}\\
            \text{and }\exists(s,t,r)\in\mathcal{I}\text{ s.t. }x_s=-1\ \&\ x_t=-\alpha_r\ \&\ x_r=-\alpha_t
        \end{matrix}\right\}.
    \end{align*}
    Moreover, the set
    \begin{align*}
        \mathfrak{U}(\mathbf{k})=\left\{\mathbf{x}\in V_\mathbf{k}(\mathbb{Z}):\begin{matrix}
            \exists i\in\{1,2,3\}\text{ s.t. }|x_i|\geq1;\\
            \exists s\in\{1,2,3\}\text{ s.t. }\forall (s,t,r)\in\mathcal{I},\ |\alpha_t-x_rx_s-x_t|<|x_t|
        \end{matrix}\right\}
    \end{align*}
    is finite, and if $\mathfrak{U}(\mathbf{k})=\emptyset$ then $\mathfrak{S}(\mathbf{k})\cup\mathfrak{T}_0(\mathbf{k})\cup\mathfrak{T}_{1}(\mathbf{k})\cup\mathfrak{T}_{-1}(\mathbf{k})$ is a fundamental domain for $\Gamma_\mathbf{k}$.
\end{theorem}
Replacing $\Gamma_\mathbf{k}$ with the mapping class group, one also has a similar result to Theorem \ref{1.1}.
\begin{theorem}\label{1.2}
    For each $\mathbf{k}\in\mathbb{C}^4$, it is decidable whether or not $\mathbf{x},\mathbf{y}\in V_\mathbf{k}(\mathbb{Z})$ are in the same $\Gamma_\mathbf{k}$-orbit. (That is, there exists an algorithm which terminates in finitely many steps and has outputs ``yes'' or ``no'' according to whether $\mathbf{x},\mathbf{y}\in V_\mathbf{k}(\mathbb{Z})$ are in the same $\Gamma_\mathbf{k}$-orbit.)
\end{theorem}
\begin{theorem}\label{1.3}
    For each $\mathbf{k}\in\mathbb{C}^4$, it is decidable whether or not $\mathbf{x},\mathbf{y}\in V_\mathbf{k}(\mathbb{Z})$ are in the same mapping class group orbit.
\end{theorem}

Subsections 1.2 and 1.3 provide a brief summary on the character varieties on smooth compact oriented surfaces, the Vieta involution maps, and the mapping class groups. In Section 2, we define a height function and introduce graphs associated with this height function, similar to what Ghosh and Sarnak have done in \cite{gs}. Additionally, we classify the local shapes of the graph near each vertex. In Section 3, we prove Theorems \ref{1.1}, \ref{1.2}, and \ref{1.3}.

\subsection{The character varieties and the mapping class groups\cite{g}}
Let $\Sigma_{g,n}$ be a smooth compact oriented surface of genus $g$ with boundary curves $\kappa_1,\dots,\kappa_n$. 
\begin{definition}[the character variety] For $\mathbf{k}=(k_1,\dots,k_n)\in\mathbb{C}^n$, the \textit{character variety} $X_\mathbf{k}$ on $\Sigma_{g,n}$ over $\mathbb{C}$ is the geometric invariant-theoretic quotient (satisfies the boundary condition $\tr\rho(\kappa_i)=k_i$ for all $i=1,\dots,n$ and $\rho\in X_\mathbf{k}$)
\begin{align*}
    X_\mathbf{k}=X_\mathbf{k}(\Sigma_{g,n})=\Hom_k(\pi_1\Sigma_{g,n},\SL_2(\mathbb{C}))\sslash\SL_2(\mathbb{C})
\end{align*}
    under the conjugation action of $\SL_2(\mathbb{C})$. For $A\subset\mathbb{C}$, let 
    \begin{align*}
        X_\mathbf{k}(\Sigma_{g,n},A)=\{\rho\in X_\mathbf{k}:\tr\rho(\gamma)\in A\ \text{for all essential curves}\ \gamma\subset\Sigma_{g,n}\},
    \end{align*}
    where the \textit{essential curve} on $\Sigma_{g,n}$ is a non-contractible simple closed curve not homotopic to $\kappa_i$'s.
\end{definition}
For given $\mathbf{k}=(k_1,\dots,k_n)\in\mathbb{C}^n$, any $\rho\in X_\mathbf{k}$ satisfies $\tr\rho(\kappa_i)=k_i$ for all $i=1,\dots,n$ as the boundary condition. Then, by the trace formula, some equations become invariants of $X_\mathbf{k}$.
\begin{definition}
    Let $(\alpha_1,\beta_1\dots,\alpha_g,\beta_g,\kappa_1,\dots,\kappa_n)$ be a sequence of generators of the fundamental group
    \begin{align*}
    \pi_1(\Sigma_{g,n},x)=\langle\alpha_1,\bar{\beta_1}\dots,\alpha_g,\bar{\beta_g},\kappa_1,\dots,\kappa_n:[\alpha_1,\bar{\beta_1}]\cdots[\alpha_g,\bar{\beta_g}]\kappa_1\cdots\kappa_n\rangle,
\end{align*}
where $\kappa_1,\dots,\kappa_n$ correspond to the boundary curves of $\Sigma_{g,n}$ and $\bar{\beta_i}$ denotes the based loop traversing $\beta_i$ in the opposite direction. If $\alpha_i$, $\beta_i$, and $\kappa_j$ are simple loops intersecting pairwise only at $x\in\Sigma_{g,n}$ and any products of generators in the sequence are simple, then we say the sequence is \textit{optimal}. 
\end{definition}

\begin{example}[one-holed torus] \label{ex:1.5} Let $k\in\mathbb{C}$. Consider an optimal sequence $(\alpha,\beta,\kappa)$ of generators of $\pi_1(\Sigma_{1,1})$, where $\kappa$ is the boundary curve of $\Sigma_{1,1}$. For any $\rho\in X_k$, by the trace formula,
\begin{align*}
    k=\tr\rho(\kappa)=\tr\rho(\alpha)^2+\tr\rho(\beta)^2+\tr\rho(\alpha\beta)^2-\tr\rho(\alpha)\tr\rho(\beta)\tr\rho(\alpha\beta)-2.
\end{align*}
Thus, writing $\rho=(x,y,z)=(\tr\rho(\alpha_1),\tr\rho(\alpha_2),\tr\rho(\alpha_1\alpha_2))$, $X_k$ is a cubic affine surface
\begin{align*}
    x^2+y^2+z^2-xyz-2-k\subset\mathbb{A}_{x,y,z}^3.
\end{align*}
\end{example}

\begin{definition}[the mapping class group]
    The \textit{mapping class group} of $\Sigma_{g,n}$ is the group of isotopy classes of orientation preserving diffeomorphisms of $\Sigma_{g,n}$ with $\partial\Sigma_{g,n}$ fixed, denoted by
    \begin{align*}
        \MCG(\Sigma_{g,n})=\pi_0(\Diff^+(\Sigma_{g,n},\partial\Sigma_{g,n})).
    \end{align*}
\end{definition}
\begin{remark}
    If $\mathbf{k}\in \mathbb{Z}^n$, then integral points in $X_\mathbf{k}$ correspond to the local systems whose monodromy elements all have integral trace. Whang showed in \cite{w} that for any $\mathbf{k}\in\mathbb{C}^n$  and $3g+n-3>0$, $X_\mathbf{k}(\Sigma_{g,n},\mathbb{Z}\setminus\{\pm2\})$ consists of at most finitely many $\MCG(\Sigma_{g,n})$-orbits.
\end{remark}

\subsection{Character variety and the mapping class group of a four-holed sphere} Let $(\gamma_1,\gamma_2,\gamma_3,\gamma_4)$ be an optimal sequence of generators for $\pi_1(\Sigma_{0,4})$, as left in Figure 1. 
   \begin{figure}[H]
\centering
   \includegraphics[width=10cm]{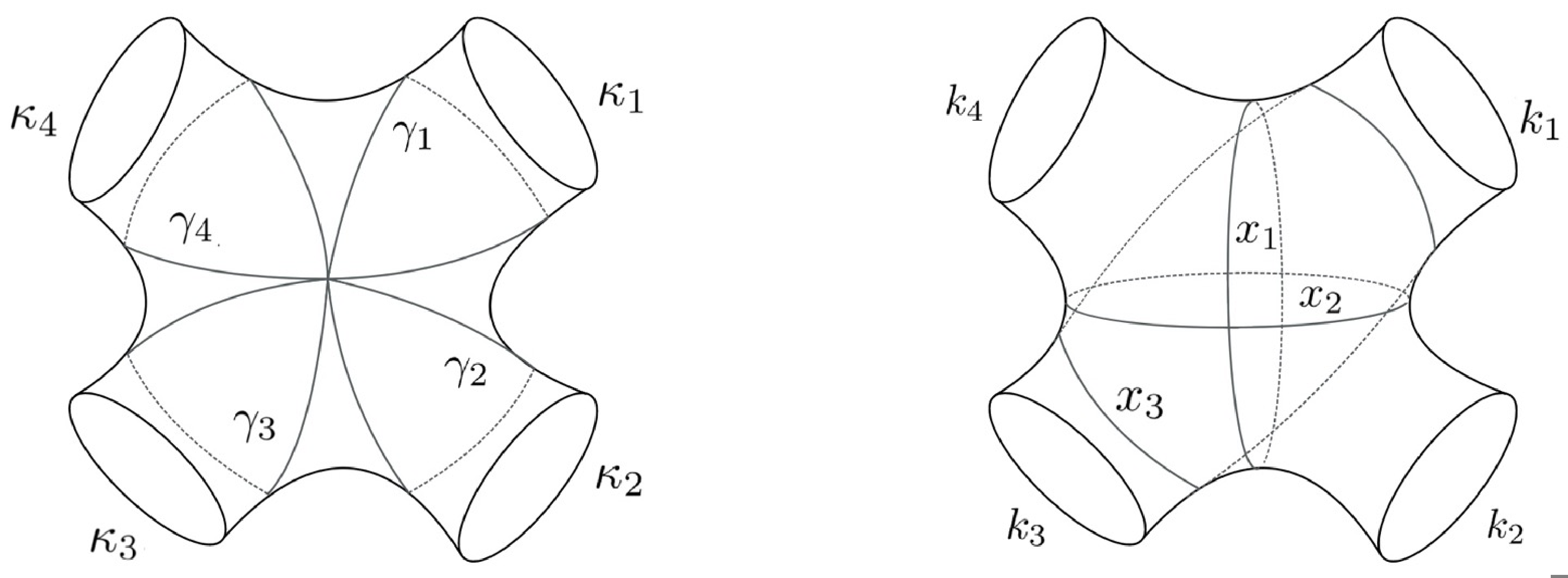}\caption{\small Curves and traces on a four-holed sphere}
\end{figure}

\noindent For each $\mathbf{k}=(k_1,k_2,k_3,k_4)\in \mathbb{C}^4$, setting $x_1=\tr\rho(\gamma_1\gamma_2),\ x_2=\tr\rho(\gamma_2\gamma_3),\  x_3=\tr\rho(\gamma_1\gamma_3)$ for some $\rho\in X_\mathbf{k}$, 
the moduli space $X_\mathbf{k}$ is a cubic affine surface 
\begin{align*}
    V_\mathbf{k}:x_1^2+x_2^2+x_3^2+x_1x_2x_3-\alpha_1x_1-\alpha_2x_2-\alpha_3x_3-\beta\subset\mathbb{A}_{x_1,x_2,x_3}^3,
\end{align*}
where
\begin{align*}
\begin{cases}
    \alpha_1=k_1k_2+k_3k_4,\\
    \alpha_2=k_1k_4+k_2k_3,\\
    \alpha_3=k_1k_3+k_2k_4,\\
    \beta=4-\sum_{i=1}^4k_i^2-\prod_{i=1}^4k_i.
\end{cases}
   \end{align*}
As mentioned before, the mapping class group $\MCG(\Sigma_{0,4})$ of a four-holed sphere acts on $V_\mathbf{k}$ via compositions of two of Vieta involution maps $\mathcal{V}_i$.

\subsection*{Acknowledgements} I thank my master thesis advisor at Seoul National University, Junho Peter Whang, for introducing me to intriguing problems and providing guidance while writing this paper. 

I was partially supported by the Samsung Science and Technology Foundation under Project Number SSTF-BA2201-03 during this work.

I would like to thank Professor Vitaly Bergelson and the (anonymous) referee for helpful comments and suggestions that improved the structure and clarity of this paper.

\section{Graphs associated with the delta map}

\subsection{Graphs associated with the delta map} Define the delta map $\Delta:\mathbb{Z}^3\rightarrow\mathbb{Z}_{\geq0}$ to be
\begin{align*}
    \Delta=\Delta(x_1,x_2,x_3)&=|x_1|+|x_2|+|x_3|.
\end{align*}
Fix $\mathbf{k}\in\mathbb{C}^4$. Define graphs $\mathcal{G}_\mathbf{x}$ associated with $\Delta$ as follows: for each $\mathbf{x}\in V_\mathbf{k}(\mathbb{Z})$, the vertices of $\mathcal{G}_\mathbf{x}$ correspond to (and labeled by) the elements of $\Gamma_\mathbf{k}\cdot\mathbf{x}$. The vertices of all graphs with the same $\Delta$-value are arbitrarily placed along the same horizontal line in $\mathbb{R}^2$, in no specific order and without overlapping; for $\mathbf{y},\mathbf{z}\in \Gamma_\mathbf{k}\cdot\mathbf{x}$ with $\Delta(\mathbf{y})<\Delta(\mathbf{z})$, the horizontal line passing through $\mathbf{y}$ is located below the one passing through $\mathbf{z}$. That is, only the values of $\Delta$ matter, not the individual coordinates. Each vertex $\mathbf{y}$ is connected with vertices $\mathcal{V}_i(\mathbf{y})$ by edges indexed by $\mathcal{V}_i$.
\begin{example}\label{ex:2.1}
    Let $k=(1,1,1,1)$. Then $(0,1,2)\in V_k(\mathbb{Z})$, and the graph corresponding to the orbit $\Gamma_{(1,1,1,1)}\cdot(0,1,2)$ is as below:
    \begin{align*}
    \begin{tikzpicture}[baseline=(2.base), node distance=7mm and 5mm, V/.style={font=\footnotesize}]
        \node[V] (1) {$(0,1,2)$};
        \node[V] (2) [below right=of 1] {$(0,1,0)$};
        \node[V] (3) [above right=of 2] {$(2,1,0)$};
        \draw (2) edge["$\mathcal{V}_3$"] (1)
              (3) edge["$\mathcal{V}_1$"] (2);
    \end{tikzpicture}
\end{align*}
\end{example}
\begin{definition}\label{def:2.2}  Let $\mathcal{G}$ be a graph associated with $\Delta$. Let $\mathbf{x}$ and $\mathbf{y}$ be vertices of $\mathcal{G}$. 
\begin{itemize}
    \item If there are $l_1,\dots,l_N\in\{1,2,3\}$ such that $\mathbf{y}=\mathcal{V}_{l_N}\cdots\mathcal{V}_{l_1}(\mathbf{x})$ and  
\begin{align*}
    \Delta\mathcal{V}_{l_{t+1}}\mathcal{V}_{l_t}\cdots\mathcal{V}_{l_1}(\mathbf{x})\geq\Delta\mathcal{V}_{l_t}\cdots\mathcal{V}_{l_1}(\mathbf{x})\geq\Delta(\mathbf{x})
\end{align*}
for all $i\in\{1,2,3\}$ and $t\in\{1,\dots,N-1\}$, then we say $\mathbf{y}$ is an \textit{ascending vertex} from $\mathbf{x}$, and $\mathbf{x}$ is a \textit{descending vertex} from $\mathbf{y}$.  
\item If $\Delta(\mathbf{x})\leq\Delta\mathcal{V}_i(\mathbf{x})$ for all $i\in\{1,2,3\}$, then we say $\mathbf{x}$ is a \textit{minimal vertex} of $\mathcal{G}$.  
\end{itemize}
\end{definition}
\begin{example}
    In Example \ref{ex:2.1}.  $(0,1,2)$ and $(2,1,0)$ are ascending vertices from $(0,1,0)$, while $(0,1,0)$ is a descending vertex from both  $(0,1,2)$ and $(2,1,0)$, and at the same time, $(0,1,0)$ is a minimal vertex of the graph.
\end{example}

\subsection{Classification of the local shapes of the graphs} Let $\mathbf{x}=(x_1,x_2,x_3)\in V_\mathbf{k}(\mathbb{Z})$. 
\begin{itemize}
    \item If $x_1=0$ then we have $\mathcal{V}_3\mathcal{V}_2(\mathbf{x})=\mathcal{V}_2\mathcal{V}_3(\mathbf{x})=(0,\alpha_2-x_2,\alpha_3-x_3)$;
    \begin{figure}[H]
    \centering
    \tikzset{
        V/.style={font=\footnotesize, align=center}, 
        every edge quotes/.style={auto, font=\footnotesize}
    }
    \begin{minipage}[b]{0.45\textwidth}
        \centering
        \begin{tikzpicture}[node distance=9mm and 6mm]
            \node[V] (1) {$(0,x_2,x_3)$};
            \node[V] (2) [below left=of 1] {$(0,\alpha_2 - x_2,x_3)$};
            \node[V] (3) [below right=of 1] {$(0,x_2,\alpha_3 - x_3)$};
            \node[V] (4) [below=1.8cm of 1] {$(0,\alpha_2 - x_2,\alpha_3 - x_3)$};
            \draw (2) edge["$\mathcal{V}_2$"] (1) (1) edge["$\mathcal{V}_3$"] (3)
                  (4) edge["$\mathcal{V}_3$"] (2) (3) edge["$\mathcal{V}_2$"] (4);
        \end{tikzpicture}
    \end{minipage}\hfill
    \begin{minipage}[b]{0.23\textwidth}
        \centering
        \begin{tikzpicture}[node distance=9mm]
            \node[V] (1) {$(0,x_2,x_3)$ \\ $=\mathcal{V}_2(0,x_2,x_3)$};
            \node[V] (3) [below=of 1] {$(0,x_2,\alpha_3 - x_3)$};
            \draw (1) edge["$\mathcal{V}_3$"] (3);
        \end{tikzpicture}
    \end{minipage}\quad
    \begin{minipage}[b]{0.23\textwidth}
        \centering
        \begin{tikzpicture}
            \node[V] (1) {$(0,x_2,x_3)$ \\ $=\mathcal{V}_2(0,x_2,x_3)$ \\ $=\mathcal{V}_3(0,x_2,x_3)$};
        \end{tikzpicture}
    \end{minipage}
    \caption{}\label{fig:0}
\end{figure}
\noindent As shown in Figure \ref{fig:0}, ascending vertices from $\mathbf{x}$, $\mathcal{V}_i(\mathbf{x})$, and $\mathcal{V}_j\mathcal{V}_i(\mathbf{x})$ with $\{i,j\}=\{2,3\}$ are $\Gamma_\mathbf{k}$-equivalent to $(0,x_2,x_3)$, $(0,\alpha_2-x_2,x_3)$, $(0,x_2,\alpha_3-x_3)$, and $(0,\alpha_2-x_2,\alpha_3-x_3)$.
\item If $x_1=1$ then we have $\mathcal{V}_2\mathcal{V}_3\mathcal{V}_2(\mathbf{x})=\mathcal{V}_3\mathcal{V}_2\mathcal{V}_3(\mathbf{x})=(1,\alpha_2-\alpha_3+x_3,\alpha_3-\alpha_2+x_2)$;
\begin{figure}[H]
\centering
\resizebox{\hsize}{!}{
\tikzset{
    lab/.style={font=\footnotesize},
    every edge quotes/.style={auto,font=\footnotesize}
}

\begin{tikzpicture}[node distance=11mm and 2mm]

\node (1) {\footnotesize$(1,x_2,x_3)$};
\node (2) [below left=of 1]
    {\footnotesize$(1,\alpha_2-x_3-x_2,x_3)$};
\node (3) [below right=of 1]
    {\footnotesize$(1,x_2,\alpha_3-x_2-x_3)$};

\node (4) [below=of 2]
    {\footnotesize$(1,\alpha_2-x_3-x_2,\alpha_3-\alpha_2+x_2)$};

\node (5) [below=of 3]
    {\footnotesize$(1,\alpha_2-\alpha_3+x_3,\alpha_3-x_2-x_3)$};

\node (7) [below=of 1] {};
\node (8) [below=of 7] {};

\node (6) [below=of 8]
    {\footnotesize$(1,\alpha_2-\alpha_3+x_3,\alpha_3-\alpha_2+x_2)$};

\path
(2) edge["$\mathcal{V}_2$"] (1)
(1) edge["$\mathcal{V}_3$"] (3)
(4) edge["$\mathcal{V}_3$"] (2)
(3) edge["$\mathcal{V}_2$"] (5)
(5) edge["$\mathcal{V}_3$"] (6)
(6) edge["$\mathcal{V}_2$"] (4);

\end{tikzpicture}

\begin{tikzpicture}[node distance=11mm and 5mm]

\node (1) {
\footnotesize$
\begin{matrix}
(1,x_2,x_3)\\
=\mathcal{V}_2(1,x_2,x_3)
\end{matrix}
$
};

\node (3) [below=of 1]
    {\footnotesize$(1,x_2,\alpha_3-x_2-x_3)$};

\node (5) [below=of 3]
    {\footnotesize$(1,\alpha_2-\alpha_3+x_3,\alpha_3-\alpha_3+x_3)$};

\path
(1) edge["$\mathcal{V}_3$"] (3)
(3) edge["$\mathcal{V}_2$"] (5);

\end{tikzpicture}

\hspace{12pt}

\begin{tikzpicture}[node distance=9mm and 6mm]

\node {
\footnotesize$
\begin{matrix}
(0,x_2,x_3)\\
=\mathcal{V}_2(0,x_2,x_3)\\
=\mathcal{V}_3(0,x_2,x_3)
\end{matrix}
$
};

\end{tikzpicture}
}
\caption{}
\label{fig:1}
\end{figure}
\noindent In Figure \ref{fig:1} above, ascending vertices from $\mathbf{x}$, $\mathcal{V}_i(\mathbf{x})$, $\mathcal{V}_j\mathcal{V}_i(\mathbf{x})$, and $\mathcal{V}_i\mathcal{V}_j\mathcal{V}_i(\mathbf{x})$ with $\{i,j\}=\{2,3\}$ are $\Gamma_\mathbf{k}$-equivalent to $(1,x_2,x_3)$, $(1,\alpha_2-x_3-x_2,x_3)$, $(1,x_2,\alpha_3-x_2-x_3)$, $(1,\alpha_2-x_3-x_2,\alpha_3-\alpha_2+x_2)$, $(1,\alpha_2-\alpha_3+x_3,\alpha_3-x_2-x_3)$, and $(1,\alpha_2-\alpha_3-x_3,\alpha_3-\alpha_2-x_2)$.

\item If $x_1=-1$ then we have $\mathcal{V}_2\mathcal{V}_3\mathcal{V}_2(\mathbf{x})=\mathcal{V}_3\mathcal{V}_2\mathcal{V}_3(\mathbf{x})=(-1,\alpha_2+\alpha_3-x_3,\alpha_3+\alpha_2-x_2)$;
\begin{figure}[H]
\centering
\resizebox{\hsize}{!}{
\tikzset{
    every edge quotes/.style={auto,font=\footnotesize}
}

\begin{tikzpicture}[node distance=11mm and 2mm]

\node (1) {\footnotesize$(-1,x_2,x_3)$};

\node (2) [below left=of 1]
    {\footnotesize$(-1,\alpha_2+x_3-x_2,x_3)$};

\node (3) [below right=of 1]
    {\footnotesize$(-1,x_2,\alpha_3+x_2-x_3)$};

\node (4) [below=of 2]
    {\footnotesize$(-1,\alpha_2+x_3-x_2,\alpha_3+\alpha_2-x_2)$};

\node (5) [below=of 3]
    {\footnotesize$(-1,\alpha_2+\alpha_3-x_3,\alpha_3+x_2-x_3)$};

\node (7) [below=of 1] {};
\node (8) [below=of 7] {};

\node (6) [below=of 8]
    {\footnotesize$(-1,\alpha_2+\alpha_3-x_3,\alpha_3+\alpha_2-x_2)$};

\path
(2) edge["$\mathcal{V}_2$"] (1)
(1) edge["$\mathcal{V}_3$"] (3)
(4) edge["$\mathcal{V}_3$"] (2)
(3) edge["$\mathcal{V}_2$"] (5)
(5) edge["$\mathcal{V}_3$"] (6)
(6) edge["$\mathcal{V}_2$"] (4);

\end{tikzpicture}

\begin{tikzpicture}[node distance=11mm and 5mm]

\node (1) {
\footnotesize$
\begin{matrix}
(-1,x_2,x_3)\\
=\mathcal{V}_2(-1,x_2,x_3)
\end{matrix}
$
};

\node (3) [below=of 1]
    {\footnotesize$(-1,x_2,\alpha_3+x_2-x_3)$};

\node (5) [below=of 3]
    {\footnotesize$(-1,\alpha_2+\alpha_3-x_3,\alpha_3+\alpha_2-x_2)$};

\path
(1) edge["$\mathcal{V}_3$"] (3)
(3) edge["$\mathcal{V}_2$"] (5);

\end{tikzpicture}

\hspace{12pt}

\begin{tikzpicture}[node distance=11mm and 5mm]

\node {
\footnotesize$
\begin{matrix}
(-1,x_2,x_3)\\
=\mathcal{V}_2(-1,x_2,x_3)\\
=\mathcal{V}_3(-1,x_2,x_3)
\end{matrix}
$
};

\end{tikzpicture}
}
\caption{}
\label{fig:-1}
\end{figure}
\noindent Ascending vertices from $\mathbf{x}$, $\mathcal{V}_i(\mathbf{x})$, $\mathcal{V}_j\mathcal{V}_i(\mathbf{x})$, and $\mathcal{V}_i\mathcal{V}_j\mathcal{V}_i(\mathbf{x})$ with $\{i,j\}=\{2,3\}$ are $\Gamma_\mathbf{k}$-equivalent to $(-1,x_2,x_3)$, $(-1,\alpha_2+x_3-x_2,x_3)$, $(-1,x_2,\alpha_3+x_2-x_3)$, $(-1,\alpha_2+x_3-x_2,\alpha_3+\alpha_2-x_2)$, $(-1,\alpha_2+\alpha_3-x_3,\alpha_3+x_2-x_3)$, and $(-1,\alpha_2+\alpha_3-x_3,\alpha_3+\alpha_2-x_2)$.
\end{itemize}
\begin{remark}\label{nf}
    If $\alpha_i-x_jx_k-x_i=x_i$, then the integer solutions $\mathbf{x}$ and $\mathcal{V}_i(\mathbf{x})$ represent the same vertex in the correspond graph. Moreover, this graph can be obtained from the graph constructed under the assumption $\alpha_i-x_jx_k-x_i\neq x_i$ by contracting the vertices $x_i$ and $\alpha_i-x_jx_k-x_i$, that is, by identifying them together with all their incident edges (see Figures \ref{fig:0}, \ref{fig:1}, and \ref{fig:-1}). Therefore, it suffices to work under assumption $\alpha_i-x_jx_k-x_i\neq x_i$, equivalently, $\mathcal{V}_i(\mathbf{x})\neq\mathbf{x}$ for all $i=1,2,3$.
\end{remark}
\begin{lemma}\label{2.1}
    If $\mathbf{y}\in V_\mathbf{k}(\mathbb{Z})$ is an ascending vertex from $\mathbf{x}=(x_1,x_2,x_3)\in V_\mathbf{k}(\mathbb{Z})$ with $|x_s|\leq1$ for some $s\in\{1,2,3\}$ in a graph $\mathcal{G}$, then $\mathbf{y}$ is $\Gamma_\mathbf{k}$-equivalent to an element in $\mathfrak{T}_0(\mathbf{k})\cup\mathfrak{T}_1(\mathbf{k})\cup\mathfrak{T}_{-1}(\mathbf{k})$, where $\mathfrak{T}_0(\mathbf{k})$, $\mathfrak{T}_1(\mathbf{k})$, and $\mathfrak{T}_{-1}(\mathbf{k})$ are the same as in Theorem \ref{1.1}.
\end{lemma}
\begin{proof}
    Without loss of generality, assume that $\mathbf{y}$ is an ascending vertex from $\mathbf{x}=(x_1,x_2,x_3)$ with $|x_1|\leq1$. As noted in Remark \ref{nf}, we may assume that $\alpha_i-x_jx_k-x_i\neq x_i$ for all $i=1,2,3$.
    \begin{itemize}
        \item Suppose $x_1=0$. Applying the same argument which was spelled out in the discussion following Figure \ref{fig:0}, we may assume that for each $t=2,3$,
        \begin{equation}\label{eq:0}
            |x_t|=\min\{|x_t|,|\alpha_t-x_t|\}.
        \end{equation}
        If $\alpha_t\neq0$ for all $t=2,3$, then $\alpha_t-x_t\neq -x_t$ and $\Delta(\mathbf{x})<\Delta\mathcal{V}_t(\mathbf{x})$. Then $(0,x_2,x_3)$ is the only one satisfying \ref{eq:0} among $\mathbf{x}$, $\mathcal{V}_i(\mathbf{x})$, and $\mathcal{V}_j\mathcal{V}_i(\mathbf{x})$ with $\{i,j\}=\{2,3\}$:
        \begin{figure}[H]
    \centering
    \begin{tikzpicture}[
        node distance = 9mm and 6mm,
        V/.style = {},
        every edge quotes/.style = {auto, font=\footnotesize}
    ]
        \begin{scope}[nodes=V]
            \node (1)   {\footnotesize{$(0,x_2,x_3)$}};
            \node (2) [above left = of 1]    {\footnotesize{$(0,\alpha_2 - x_2,x_3)$}};
            \node (3) [above right = of 1]   {\footnotesize{$(0,x_2,\alpha_3 - x_3)$}};
            \node (5) [above = of 1]         {};
            \node (4) [above = of 5]         {\footnotesize{$(0,\alpha_2 - x_2,\alpha_3 - x_3)$}};
        \end{scope}
        \draw 
            (1) edge["$\mathcal{V}_2$"] (2)
            (3) edge["$\mathcal{V}_3$"] (1)
            (2) edge["$\mathcal{V}_3$"] (4)
            (4) edge["$\mathcal{V}_2$"] (3);
    \end{tikzpicture}
\end{figure}
\noindent Thus,
\begin{align*}
    \mathbf{x}\in\left\{\mathbf{y}\in V_\mathbf{k}(\mathbb{Z}):
            \begin{matrix}
                \exists s\in\{1,2,3\}\text{ s.t. }y_s=0\\
            \&\ \forall t\neq s,\ |y_t|=\min\{|y_t|,|\alpha_t-y_t|\}\ \&\ \alpha_t\neq0
            \end{matrix}\right\}.
\end{align*}
If $\alpha_2=0$ and $\alpha_3\neq0$, then $\alpha_2-x_2= -x_2$, $\alpha_3-x_3\neq -x_3$, and thus $\Delta(\mathbf{x})=\Delta\mathcal{V}_2(\mathbf{x})$, $\Delta(\mathbf{x})<\Delta\mathcal{V}_3(\mathbf{x})$:
        \begin{figure}[H]
    \centering
    \begin{tikzpicture}[
        node distance = 10mm and 10mm,
        V/.style = {},
        every edge quotes/.style = {auto, font=\footnotesize}
    ]
        \begin{scope}[nodes=V]
            \node (1)   {\footnotesize{$(0,x_2,x_3)$}};
            \node (2) [left = of 1]    {\footnotesize{$(0,-x_2,x_3)$}};
            \node (3) [above = of 1]   {\footnotesize{$(0,x_2,\alpha_3 - x_3)$}};
            \node (4) [above = of 2]   {\footnotesize{$(0,-x_2,\alpha_3 - x_3)$}};
        \end{scope}
        
        \draw 
            (1) edge["$\mathcal{V}_2$"] (2)
            (3) edge["$\mathcal{V}_3$"] (1)
            (2) edge["$\mathcal{V}_3$"] (4)
            (4) edge["$\mathcal{V}_2$"] (3);
    \end{tikzpicture}
\end{figure}
\noindent Since $(0,\pm x_2,x_3)$ are $\Gamma_\mathbf{k}$-equivalent, we may assume $|x_2|=x_2$, which guarantees that the minimal vertex within $\mathfrak{T}_0(\mathbf{k})$ is unique on the corresponding graph when $\mathfrak{U}(\mathbf{k})=\emptyset$. Then 
\begin{align*}
    \mathbf{x}\in\left\{\mathbf{y}\in V_\mathbf{k}(\mathbb{Z}):
            \begin{matrix}
                \exists (s,t,r)\in\mathcal{I}\text{ s.t. }y_s=0\ \&\ |y_t|=y_t\\
            \&\ |y_r|=\min\{|y_r|,|\alpha_r-y_r|\}\ \&\ \alpha_t=0\ \&\ \alpha_r\neq0
            \end{matrix}\right\}.
\end{align*}
        If $\alpha_2=\alpha_3=0$, then $\alpha_t-x_t= -x_t$ and $\Delta(\mathbf{x})=\Delta\mathcal{V}_t(\mathbf{x})$ for $t=2,3$. Since $(0,\pm x_2,\pm x_3)$  are $\Gamma_\mathbf{k}$-equivalent, we may assume $|x_t|=x_t$ for $t=2,3$. Thus,
        \begin{align*}
    \mathbf{x}\in\left\{\mathbf{y}\in V_\mathbf{k}(\mathbb{Z}):
            \exists s\in\{1,2,3\}\text{ s.t. } y_s=0\ \&\ \forall t\neq s,\ |y_t|=y_t\ \&\ \alpha_t=0\right\}.
\end{align*}
\item Suppose $x_1=1$. We may assume that for each $(t,r)\in\{(2,3),(3,2)\}$,
        \begin{equation}\label{eq:1}
            1\leq|x_t|=\min\{|x_t|,|\alpha_t-x_r-x_t|,|\alpha_t-\alpha_r-x_r|\}.
        \end{equation}
    If $\mathcal{V}_1(\mathbf{z})=(0,z_2,z_3)$, where $\mathbf{z}=(z_1,z_2,z_3)$ is one of $\mathbf{x}$, $\mathcal{V}_i(\mathbf{x})$, $\mathcal{V}_j\mathcal{V}_i(\mathbf{x})$, and $\mathcal{V}_i\mathcal{V}_j\mathcal{V}_i(\mathbf{x})$ with $\{i,j\}=\{2,3\}$, then $\mathbf{x}$ is an ascending vertex from $\mathbf{w}=(0,w_2,w_3)$; by the previous argument with Figure \ref{fig:1}, $\mathbf{x}$ is $\Gamma_\mathbf{k}$-equivalent to one in $\mathfrak{T}_0(\mathbf{k})$. Now, we may assume that for $(t,r)\in\{(2,3),(3,2)\}$,
\begin{equation}\label{eq:2}
    1\leq\min\left\{\begin{matrix}
                |\alpha_1-x_tx_r-1|,|\alpha_1-(\alpha_t-x_r-x_t)x_r-1|,\\
                |\alpha_1-(\alpha_t-x_r-x_t)(\alpha_r-\alpha_t+x_t)-1|,\\
                |\alpha_1-(\alpha_t-\alpha_r+x_r)(\alpha_r-\alpha_t+x_t)-1|
            \end{matrix}\right\},
\end{equation}
i.e.,\ all of $\mathbf{x}$, $\mathcal{V}_i(\mathbf{x})$, $\mathcal{V}_j\mathcal{V}_i(\mathbf{x})$, and $\mathcal{V}_i\mathcal{V}_j\mathcal{V}_i(\mathbf{x})$ with $\{i,j\}=\{2,3\}$ have no descending vertex of the form $\mathbf{w}=(0,w_2,w_3)$ (Note that inequality \ref{eq:2} means that
\begin{align*}
    \Delta(\mathbf{x})\leq\Delta\mathcal{V}_1(\mathbf{x}),\Delta\mathcal{V}_1\mathcal{V}_t,\Delta\mathcal{V}_1\mathcal{V}_t\mathcal{V}_r(\mathbf{x}),\Delta\mathcal{V}_1\mathcal{V}_t\mathcal{V}_r\mathcal{V}_t(\mathbf{x}),
\end{align*}
so that the first coordinate cannot become zero under the action of $\mathcal{V}_1$; see the vertices--especially the second and third coordinates--in Figure \ref{fig:1}. Indeed, the inequality \ref{eq:2} is introduced precisely to prevent $\mathbf{x}$ from having a descending vertex of the form $(0,w_2,w_3)$, since we are collecting minimal vertices $\mathbf{x}$ satisfying $\min_{s=1,2,3}|x_s|=1=x_i$). If $x_t\neq\alpha_r$ for $(t,r)\in\{(2,3),(3,2)\}$, then $\alpha_r-x_t-x_r\neq-x_r$ and $\Delta(\mathbf{x})<\Delta\mathcal{V}_r(\mathbf{x})$. Then $(1,x_2,x_3)$ is the only one satisfying \ref{eq:1} among $\mathbf{x}$, $\mathcal{V}_i(\mathbf{x})$, $\mathcal{V}_j\mathcal{V}_i(\mathbf{x})$, and $\mathcal{V}_i\mathcal{V}_j\mathcal{V}_i(\mathbf{x})$ with $\{i,j\}=\{2,3\}$: 
\begin{figure}[H]
\centering

\tikzset{
    every edge quotes/.style={auto,font=\footnotesize}
}

\begin{tikzpicture}[node distance=11mm and 5mm]

\node (1) {\footnotesize$(1,x_2,x_3)$};

\node (2) [above left=of 1]
    {\footnotesize$(1,\alpha_2-x_3-x_2,x_3)$};

\node (3) [above right=of 1]
    {\footnotesize$(1,x_2,\alpha_3-x_2-x_3)$};

\node (4) [above=of 2]
    {\footnotesize$(1,\alpha_2-x_3-x_2,\alpha_3-\alpha_2+x_2)$};

\node (5) [above=of 3]
    {\footnotesize$(1,\alpha_2-\alpha_3+x_3,\alpha_3-x_2-x_3)$};

\node (7) [above=of 1] {};
\node (8) [above=of 7] {};

\node (6) [above=of 8]
    {\footnotesize$(1,\alpha_2-\alpha_3+x_3,\alpha_3-\alpha_2+x_2)$};

\path
(1) edge["$\mathcal{V}_2$"] (2)
(3) edge["$\mathcal{V}_3$"] (1)
(2) edge["$\mathcal{V}_3$"] (4)
(5) edge["$\mathcal{V}_2$"] (3)
(6) edge["$\mathcal{V}_3$"] (5)
(4) edge["$\mathcal{V}_2$"] (6);

\end{tikzpicture}
\end{figure}
  \noindent  So, we have
    \begin{align*}
    \mathbf{x}\in\left\{\mathbf{y}\in V_\mathbf{k}(\mathbb{Z}):\begin{matrix}\exists s\in\{1,2,3\}\text{ s.t. }\forall(s,t,r)\in\mathcal{I},\ y_t\neq\alpha_r\ \\
            \&\ y_s=1\leq\min\left\{\begin{matrix}
                |\alpha_s-y_ty_r-1|,|\alpha_s-(\alpha_t-y_r-y_t)y_r-1|,\\
                |\alpha_s-(\alpha_t-y_r-y_t)(\alpha_r-\alpha_t+y_t)-1|,\\
                |\alpha_s-(\alpha_t-\alpha_r+y_r)(\alpha_r-\alpha_t+y_t)-1|
            \end{matrix}\right\}\\
            \&\ 1\leq|y_t|=\min\{|y_t|,|\alpha_t-y_r-y_t|,|\alpha_t-\alpha_r-y_r|\}
        \end{matrix}\right\}.
\end{align*}
If $x_3=\alpha_2$ and $x_2\neq\alpha_3$, then $\alpha_2-x_3-x_2= -x_2$, $\alpha_3-x_2-x_3\neq -x_3$, and thus $\Delta(\mathbf{x})=\Delta\mathcal{V}_2(\mathbf{x})$, $\Delta(\mathbf{x})<\Delta\mathcal{V}_3(\mathbf{x})$:
\begin{figure}[H]
\centering

\tikzset{
    every edge quotes/.style={auto,font=\footnotesize}
}

\begin{tikzpicture}[node distance=13mm and 5mm]

\node (1) {\footnotesize$(1,x_2,\alpha_2)$};

\node (2) [left=of 1]
    {\footnotesize$(1,-x_2,\alpha_2)$};

\node (3) [above right=of 1]
    {\footnotesize$(1,x_2,\alpha_3-\alpha_2-x_2)$};

\node (4) [above=of 2]
    {\footnotesize$(1,-x_2,\alpha_3-\alpha_2+x_2)$};

\node (5) [above=of 3]
    {\footnotesize$(1,2\alpha_2-\alpha_3,\alpha_3-\alpha_2-x_2)$};

\node (7) [above=of 1] {};
\node (8) [above=of 7] {};

\node (6) [above=of 8]
    {\footnotesize$(1,2\alpha_2-\alpha_3,\alpha_3-\alpha_2+x_2)$};

\path
(1) edge["$\mathcal{V}_2$"] (2)
(3) edge["$\mathcal{V}_3$"] (1)
(2) edge["$\mathcal{V}_3$"] (4)
(5) edge["$\mathcal{V}_2$"] (3)
(6) edge["$\mathcal{V}_3$"] (5)
(4) edge["$\mathcal{V}_2$"] (6);

\end{tikzpicture}
\end{figure}
    \noindent Since $(1,\pm x_2,x_3)$ are $\Gamma_\mathbf{k}$-equivalent, we may assume $|x_2|=x_2$, ensuring that the minimal vertex within $\mathfrak{T}_1(\mathbf{k})$ is unique on the corresponding graph when $\mathfrak{U}(\mathbf{k})=\emptyset$. Then 
\begin{align*}
    \mathbf{x}\in\left\{\mathbf{y}\in V_\mathbf{k}(\mathbb{Z}):\begin{matrix}\exists (s,t,r)\in\mathcal{I}\text{ s.t. }y_t\neq\alpha_r\ \&\  y_r=\alpha_t\\
            \&\ y_s=1\leq\min\left\{\begin{matrix}
                |\alpha_s\pm\alpha_ty_t-1|,|\alpha_s- y_t(y_t\pm(\alpha_r-\alpha_t))-1|,\\
                |\alpha_s\pm (2\alpha_t-\alpha_r)(\alpha_r-\alpha_t\pm y_t)-1|
            \end{matrix}\right\}\\
            \&\ 1\leq|\alpha_t|=\min\{|\alpha_t|,|\alpha_r-\alpha_t\pm y_t|\}\\
            \&\ 1\leq y_t=\min\{|y_t|,|2\alpha_t-\alpha_r|\}
        \end{matrix}\right\}.
\end{align*}
If $x_t=\alpha_r$ for all $(t,r)\in\{(2,3),(3,2)\}$, then $\alpha_t-x_r-x_t= -x_t$ and $\Delta(\mathbf{x})=\Delta\mathcal{V}_t(\mathbf{x})$. Since $(0,\pm x_2,\pm x_3)$  are $\Gamma_\mathbf{k}$-equivalent, we may assume $|x_t|=x_t$ for all $t\in\{2,3\}$. Thus,
        \begin{align*}
    \mathbf{x}\in\left\{\mathbf{y}\in V_\mathbf{k}(\mathbb{Z}):
            \exists s\in\{1,2,3\}\text{ s.t. } y_s=1\ \&\ \forall (s,t,r)\in\mathcal{I},\ |y_t|=y_t=\alpha_r\right\}.
\end{align*}
\item Suppose $x_1=-1$. Similarly to the case $x_1=1$, assume 
that for $(t,r)\in\{(2,3),(3,2)\}$,
\begin{align*}
    1\leq\min\left\{\begin{matrix}
                |\alpha_1-x_tx_r+1|,|\alpha_1-(\alpha_t-x_r-x_t)x_r+1|,\\
                |\alpha_1-(\alpha_t-x_r-x_t)(\alpha_r-\alpha_t+x_t)+1|,\\
                |\alpha_1-(\alpha_t-\alpha_r+x_r)(\alpha_r-\alpha_t+x_t)+1|
            \end{matrix}\right\};
\end{align*}
otherwise, $\mathbf{x}\in\mathfrak{T}_0(\mathbf{k})$. If there exists $s\in\{1,2,3\}$ such that $x_s=-1$ and $\alpha_s=x_tx_r$ for $(s,t,r)\in\mathcal{I}$, then $\alpha_s-x_tx_r+1=1$ and $\mathbf{x}$ is $\Gamma_\mathbf{k}$-equivalent to one in $\mathfrak{T}_1(\mathbf{k})$. Therefore, to avoid reducing to the case $x_i=1$ and to ensure the uniqueness of the minimal vertex within $\mathfrak{T}_1(\mathbf{k})\cup\mathfrak{T}_{-1}(\mathbf{k})$ on the corresponding graph when $\mathfrak{U}(\mathbf{k})=\emptyset$, we may assume that for all $s\in\{1,2,3\}$ with $x_s=-1$, $\alpha_s\neq x_tx_r$ for $(s,t,r)\in\mathcal{I}$. Applying the same argument as in Figure \ref{fig:-1}, we get $\mathbf{x}\in\mathfrak{T}_{-1}(\mathbf{k})$.
    \end{itemize}
\end{proof}
\begin{remark}
    Each vertex in $\mathfrak{T}_0(\mathbf{k})\cup\mathfrak{T}_1(\mathbf{k})\cup\mathfrak{T}_{-1}(\mathbf{k})$ is minimal; moreover, as shown in the proof of Lemma \ref{2.1}, for each $l=0,\pm1$, none of the vertices in $\mathfrak{T}_l(\mathbf{k})$ are ascending or descending vertices from each other. Therefore, $\mathbf{x},\mathbf{y}\in \mathfrak{T}_0(\mathbf{k})\cup\mathfrak{T}_1(\mathbf{k})\cup\mathfrak{T}_{-1}(\mathbf{k})$ are not $\Gamma_\mathbf{k}$-equivalent unless they have a common ascending vertex $\mathbf{z}=(z_1,z_2,z_3)$ with $|z_s|\geq1$.
\end{remark}
\begin{lemma}\label{2.2}
    Let $\mathbf{x}=(x_1,x_2,x_3)\in V_\mathbf{k}(\mathbb{Z})$, and let $\mathcal{G}_\mathbf{x}$ the graph corresponding to its orbit $\Gamma_\mathbf{k}\cdot\mathbf{x}$. If $|x_s|\geq3$ and $|x_tx_r|>3|\alpha_s|$ for all $(s,t,r)\in\mathcal{I}$ and if $\mathbf{x}$ is a minimal vertex of $\mathcal{G}_\mathbf{x}$, then $\mathcal{G}_\mathbf{x}$ has no other minimal vertices.
\end{lemma}
\begin{proof}
    Let $\mathbf{x}\in V_\mathbf{k}(\mathbb{Z})$ with $|x_s|\geq3$ and $|x_tx_r|>3|\alpha_s|$ for all $(s,t,r)\in\mathcal{I}$. 
     It suffices to assume that $\alpha_i-x_jx_k-x_i\neq x_i$ for all $i=1,2,3$, as stated in Remark \ref{nf}. Then for each $(s,t,r)\in\mathcal{I}$, $\alpha_s-x_tx_r-x_s\neq-x_s$ and $\Delta\mathcal{V}_s(\mathbf{x})\neq\Delta(\mathbf{x})$. If $|\alpha_s-x_tx_r-x_s|\leq|x_s|$ for some $(s,t,r)\in\mathcal{I}$, then 
    \begin{align*}
        |x_tx_r|&\leq\begin{cases}
            \max\{|\alpha_s-x_s|,|\alpha_s-2x_s|\}&\text{if $x_s(\alpha_s-x_tx_r-x_s)\geq0$}\\
            \max\{|\alpha_s|,|\alpha_s-x_s|\}&\text{if $x_s(\alpha_s-x_tx_r-x_s)<0$} 
        \end{cases}\\
        &\leq|\alpha_s|+2|x_s|<\frac{1}{3}|x_tx_r|+2|x_s|,
    \end{align*}
    that is, $|x_t|,|x_r|<|x_s|$, which implies that $|x_t|<|\alpha_t-x_rx_s-x_t|$ and $|x_r|<|\alpha_r-x_sx_t-x_r|$;
    \begin{figure}[H]
\centering
\resizebox{\hsize}{!}{
\tikzset{
    every edge quotes/.style={auto,font=\footnotesize}
}

\begin{tikzpicture}[node distance=5mm and 8mm]

\node (1) {\footnotesize$(x_1,x_2,x_3)$};

\node (2) [above left=of 1]
{\footnotesize$(\alpha_1-x_2x_3-x_1,x_2,x_3)$};

\node (3) [above=of 1]
{\footnotesize$(x_1,\alpha_2-x_3x_1-x_2,x_3)$};

\node (4) [above right=of 1]
{\footnotesize$(x_1,x_2,\alpha_3-x_1x_2-x_3)$};

\path
(1) edge["$\mathcal{V}_1$"] (2)
(3) edge["$\mathcal{V}_2$"] (1)
(4) edge["$\mathcal{V}_3$"] (1);

\end{tikzpicture}

\begin{tikzpicture}[node distance=5mm and 1mm]

\node (1) {\footnotesize$(x_1,x_2,x_3)$};

\node (2) [below=of 1]
{\footnotesize$(\alpha_1-x_2x_3-x_1,x_2,x_3)$};

\node (3) [above=of 1]
{\footnotesize$(x_1,\alpha_2-x_3x_1-x_2,x_3)$};

\node (4) [right=of 3]
{\footnotesize$(x_1,x_2,\alpha_3-x_1x_2-x_3)$};

\path
(1) edge["$\mathcal{V}_1$"] (2)
(1) edge["$\mathcal{V}_2$"] (3)
(4) edge["$\mathcal{V}_3$"] (1);

\end{tikzpicture}
}
\end{figure}
    \noindent Therefore, if $\mathbf{x}$ is minimal, then all the vertices $\mathcal{V}_1(\mathbf{x})$, $\mathcal{V}_2(\mathbf{x})$, and $\mathcal{V}_3(\mathbf{x})$ connected with $\mathbf{x}$ on $\mathcal{G}_\mathbf{x}$ are located above $\mathbf{x}$ on the plane as follows: 
    \begin{align*}
    \begin{tikzpicture}[
node distance = 5mm and 8mm,
     V/.style = {},
every edge quotes/.style = {auto, font=\footnotesize}
                    ]
    \begin{scope}[nodes=V]
\node (1)   {\footnotesize{$\mathbf{x}$}};
\node (2) [above left = of 1]    {\footnotesize{$\mathcal{V}_1(\mathbf{x})$}};
\node (3) [above = of 1]          {\footnotesize{$\mathcal{V}_2(\mathbf{x})$}};
\node (4) [above right = of 1]    {\footnotesize{$\mathcal{V}_3(\mathbf{x})$}};
    \end{scope}
\draw   (1)  edge["\footnotesize{$\mathcal{V}_1$}"] (2)
        (3)  edge["\footnotesize{$\mathcal{V}_2$}"] (1)
        (4)  edge["\footnotesize{$\mathcal{V}_3$}"] (1);
    \end{tikzpicture}
    \end{align*}  
Note that by Definition \ref{def:2.2}, any ascending vertex $\mathbf{y}=(y_1,y_2,y_3)$ from $\mathbf{x}$ has $|x_s|\leq|y_s|$ so that $|\alpha_s|<|x_tx_r|\leq|y_ty_r|$ for all $(s,t,r)\in\mathcal{I}$. Then $\mathcal{G}_\mathbf{x}$ is illustrated as follows:
\begin{align*}
    \begin{tikzpicture}[
node distance = 5mm and 8mm,
     V/.style = {},
every edge quotes/.style = {auto, font=\footnotesize}
                    ]
    \begin{scope}[nodes=V]
\node (1)   {\footnotesize{$\mathbf{x}$}};
\node (2) [above left = of 1]    {\footnotesize{$\mathcal{V}_1(\mathbf{x})$}};
\node (3) [above = of 1]          {\footnotesize{$\mathcal{V}_2(\mathbf{x})$}};
\node (4) [above right = of 1]    {\footnotesize{$\mathcal{V}_3(\mathbf{x})$}};
\node (5) [above left = of 2]    {\footnotesize{$\mathcal{V}_3\mathcal{V}_1(\mathbf{x})$}};
\node (6) [left= of 5]          {\footnotesize{$\mathcal{V}_2\mathcal{V}_1(\mathbf{x})$}};
\node (7) [right = of 5]    {\footnotesize{$\mathcal{V}_3\mathcal{V}_2(\mathbf{x})$}};
\node (8) [above = of 3]          {\footnotesize{$\mathcal{V}_1\mathcal{V}_2(\mathbf{x})$}};
\node (9) [right = of 8]    {\footnotesize{$\mathcal{V}_1\mathcal{V}_3(\mathbf{x})$}};
\node (10) [above right= of 4]          {\footnotesize{$\mathcal{V}_2\mathcal{V}_3(\mathbf{x})$}};
\node (11) [above = of 5]          {\footnotesize{$\cdots$}};
\node (12) [above = of 6]          {\footnotesize{$\cdots$}};
\node (13) [above = of 7]          {\footnotesize{$\cdots$}};
\node (14) [above = of 8]          {\footnotesize{$\cdots$}};
\node (15) [above = of 9]          {\footnotesize{$\cdots$}};
\node (16) [above = of 10]          {\footnotesize{$\cdots$}};
    \end{scope}
\draw   (1)  edge["\footnotesize{$\mathcal{V}_1$}"] (2)
        (3)  edge["\footnotesize{$\mathcal{V}_2$}"] (1)
        (4)  edge["\footnotesize{$\mathcal{V}_3$}"] (1)
        (5)  edge["\footnotesize{$\mathcal{V}_3$}"] (2)
        (2)  edge["\footnotesize{$\mathcal{V}_2$}"] (6)
        (3)  edge["\footnotesize{$\mathcal{V}_3$}"] (7)
        (8)  edge["\footnotesize{$\mathcal{V}_1$}"] (3)
        (4)  edge["\footnotesize{$\mathcal{V}_1$}"] (9)
        (10)  edge["\footnotesize{$\mathcal{V}_2$}"] (4)
        (11)  edge[""] (5)
        (12)  edge[""] (6)
        (13)  edge[""] (7)
        (14)  edge[""] (8)
        (15)  edge[""] (9)
        (16)  edge[""] (10);
    \end{tikzpicture}
    \end{align*}
    Therefore, $\mathbf{x}$ is the unique minimal vertex of $\mathcal{G}_\mathbf{x}$.
\end{proof}
Now, let $\mathfrak{U}(\mathbf{k})$ be the same as in Theorem \ref{1.1}, and consider
\begin{align*}
    \mathfrak{U}'(\mathbf{k})&=\{(x_1,x_2,x_3)\in \mathfrak{U}(\mathbf{k}):\exists (s,t,r)\in\mathcal{I}\text{ s.t. } |x_s|=1\text{ or } |x_tx_r|\leq3|\alpha_s|\}.
\end{align*}
From the proof of Lemma \ref{2.2}, we have
\begin{align*}
    \mathfrak{U}(\mathbf{k})\setminus\mathfrak{U}'(\mathbf{k})&=\{(x_1,x_2,x_3)\in \mathfrak{U}(\mathbf{k}):\forall (s,t,r)\in\mathcal{I},\ |x_s|\geq2\ \&\ |x_tx_r|>3|\alpha_s|\}\\
    &=\left\{(x_1,x_2,x_3)\in \mathfrak{U}(\mathbf{k}):\begin{matrix}
        \forall(i,j,k)\in\mathcal{I},\ |x_jx_k|>3|\alpha_i|;\\
        \exists s\in\{1,2,3\}\text{ s.t. } |x_s|=2\ \&\ \forall t\neq s,\ |x_t|\geq3
    \end{matrix}\right\}.
\end{align*}
\begin{lemma}\label{2.3}
    The set $\mathfrak{U}(\mathbf{k})$ is finite.
\end{lemma}
\begin{proof}
    Let $\mathfrak{U}'(\mathbf{k})$ and $\mathfrak{U}(\mathbf{k})\setminus\mathfrak{U}'(\mathbf{k})$ be the same as above.
    \begin{itemize}
        \item $\mathfrak{U}'(\mathbf{k})$ is finite: If $|x_s|=1$ for some $(s,t,r)\in\mathcal{I}$, then $V_\mathbf{k}$ becomes 
        \begin{align*}
            \left(x_t\pm\frac{1}{2}x_r\right)^2+\frac{3}{4}x_r^2=-1+\alpha_tx_t+\alpha_rx_r+\beta\pm\alpha_s,
        \end{align*}
        so that $V_\mathbf{k}$ has at most finitely many integral solutions. If $|x_tx_r|\leq3|\alpha_s|$, then for each chosen $(x_t,x_r)$ from at most finitely many candidates, thus regarding $x_t,x_r$ as constants, $V_\mathbf{k}$ becomes a quadratic form in terms of $x_s$, implying $V_\mathbf{k}$ has at most two integral solutions for given $(x_t,x_r)$. Hence, $\mathfrak{U}'(\mathbf{k})$ is finite.
        \item $\mathfrak{U}(\mathbf{k})\setminus\mathfrak{U}'(\mathbf{k})$ is finite: Let $\mathbf{x}\in\mathfrak{U}(\mathbf{k})\setminus\mathfrak{U}'(\mathbf{k})$. For all $(s,t,r)\in\mathcal{I}$, $|x_tx_r|>3|\alpha_s|$, then $\alpha_s-x_tx_r-x_s\neq-x_s$ and $\Delta\mathcal{V}_s(\mathbf{x})\neq\Delta(\mathbf{x})$. Without loss of generality, assume $|x_1|=2$ and $|x_2|,|x_3|\geq3$. If $|\alpha_1-x_2x_3-x_1|<|x_1|$, then
        \begin{align*}
            |x_2x_3|<|\alpha_1|+2|x_1|<\frac{1}{3}|x_2x_3|+4,
        \end{align*}
        that is, $|x_2x_3|<6$, which is a contradiction. Thus, $|x_1|<|\alpha_1-x_2x_3-x_1|$. As $\mathfrak{U}(\mathbf{k})$ defined, $|\alpha_2-x_3x_1-x_2|<|x_2|$ and $|\alpha_3-x_1x_2-x_3|<|x_3|$.
        \begin{figure}[H]
        \centering
    \begin{tikzpicture}[
    node distance = 4mm and 1mm,
    V/.style = {},
    every edge quotes/.style = {auto, font=\footnotesize}
]
    \begin{scope}[nodes=V]
        \node (1)   {\footnotesize{$(x_1,x_2,x_3)$}};
        \node (2) [above = of 1]    {\footnotesize{$(\alpha_1 - x_2 x_3 - x_1, x_2, x_3)$}};
        \node (3) [below left = of 1] {\footnotesize{$(x_1, \alpha_2 - x_3 x_1 - x_2, x_3)$}};
        \node (4) [below right = of 1] {\footnotesize{$(x_1, x_2, \alpha_3 - x_1 x_2 - x_3)$}};
    \end{scope}
    \draw   (2) edge["$\mathcal{V}_1$"] (1)
            (3) edge["$\mathcal{V}_2$"] (1)
            (1) edge["$\mathcal{V}_3$"] (4);
\end{tikzpicture}
    \end{figure}
    \noindent Since $(2,x_2,x_3)$ or $(-2,x_2,x_3)$ is solution for $V_\mathbf{k}$, we have
    \begin{equation}\label{eq:2b}
        4+x_2^2+x_3^2\pm2(x_2x_3-\alpha_1)-\alpha_2x_2-\alpha_3x_3-\beta=0.
    \end{equation}
    If the equation \ref{eq:2b} describes a double line, i.e.,
    \begin{align*}
       4+x_2^2+x_3^2\pm2(x_2x_3-\alpha_1)-\alpha_2x_2-\alpha_3x_3-\beta=\left(x_2\pm x_3-\frac{\alpha_2}{2}\right)^2,
    \end{align*}
    then there are only two cases as in Figure \ref{fig:dl}.
\begin{figure}[h]
     \centering
     \begin{subfigure}[h]{0.35\textwidth}
         \centering
         \resizebox{40mm}{40mm}{%
        \begin{tikzpicture}
\begin{axis}
\addplot[color=red,thick]{-x-2};
\end{axis}
\end{tikzpicture}}
         \caption{\small$x_1=2$\normalsize}
         \label{fig:12}
     \end{subfigure}
     \begin{subfigure}[h]{0.35\textwidth}
         \centering
         \resizebox{40mm}{40mm}{%
        \begin{tikzpicture}
\begin{axis}
\addplot[color=red,thick]{x-2};
\end{axis}
\end{tikzpicture}}
         \caption{\small$x_1=-2$\normalsize}
         \label{fig:1-2}
     \end{subfigure}
        \caption{\small The equation \ref{eq:2} describes a double line\normalsize}
        \label{fig:dl}
\end{figure}
\noindent Note that $\mathcal{V}_i$ sends $x_i$ to $\alpha_i-x_jx_k-x_i$ and fixes the other two $x_j$ and $x_k$. Thus, if $\mathbf{x}\in\mathfrak{U}(\mathbf{k})\setminus\mathfrak{U}'(\mathbf{k})$ and it is a point on one of the graphs in Figure \ref{fig:dl}, then the line $\left(x_2\pm x_3-\frac{\alpha_2}{2}\right)^2$ has two distinct intersection points with each of lines $x=x_2$ and $y=x_3$ because $|\alpha_2-x_3x_1-x_2|<|x_2|$ and $|\alpha_3-x_1x_2-x_3|<|x_3|$. However, it is impossible. Thus, the equation \ref{eq:2} does not represent a double line, i.e.,
    \begin{align*}
       \left(x_2\pm x_3-\frac{\alpha_2}{2}\right)^2&\neq4+x_2^2+x_3^2\pm2(x_2x_3-\alpha_1)-\alpha_2x_2-\alpha_3x_3-\beta\\
       &=\left(x_2\pm x_3-\frac{\alpha_2}{2}\right)^2+(\alpha_2-\alpha_3)x_3+\gamma,\quad(\gamma\in\mathbb{Q})
    \end{align*}
    and thus there are only six cases as illustrated by Figure \ref{fig:ndl}.
    \begin{figure}[h]
     \centering
     \hspace{18pt}\includegraphics[width=14cm]{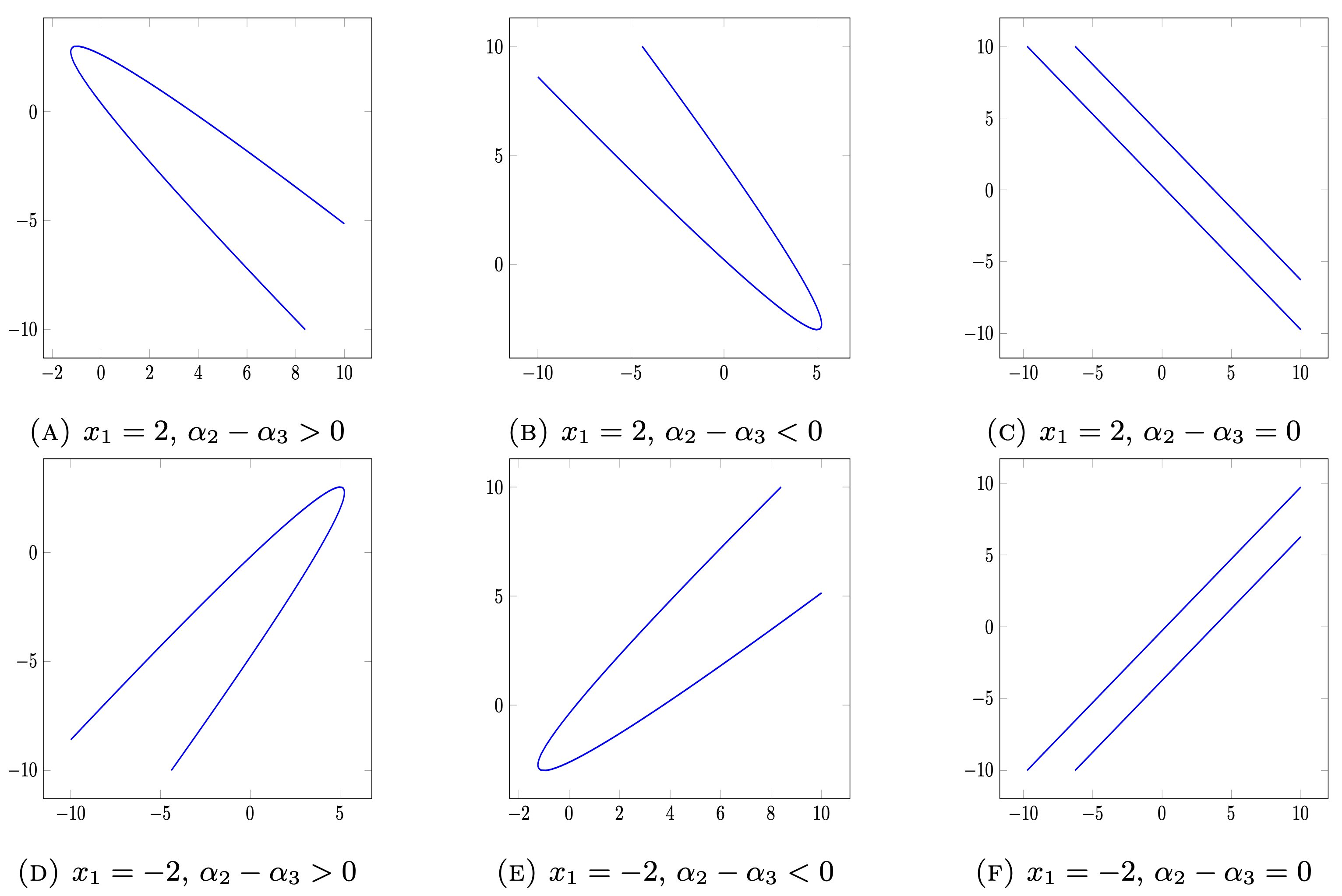}
        \caption{\small The equation \ref{eq:2} does not present a double line\normalsize}
        \label{fig:ndl}
\end{figure}
  \noindent  In all cases, there is $M>0$ such that if $|x_2|>M$ or $|x_3|>M$, then
\begin{align*}
    |\alpha_2-x_3x_1-x_2|\leq|x_2|\ \ \text{and}\ \  |\alpha_3-x_1x_2-x_3|\leq|x_3|
\end{align*}
cannot hold simultaneously (see Figure \ref{b}).
    \begin{figure}[H]
\centering
\includegraphics[width=4cm]{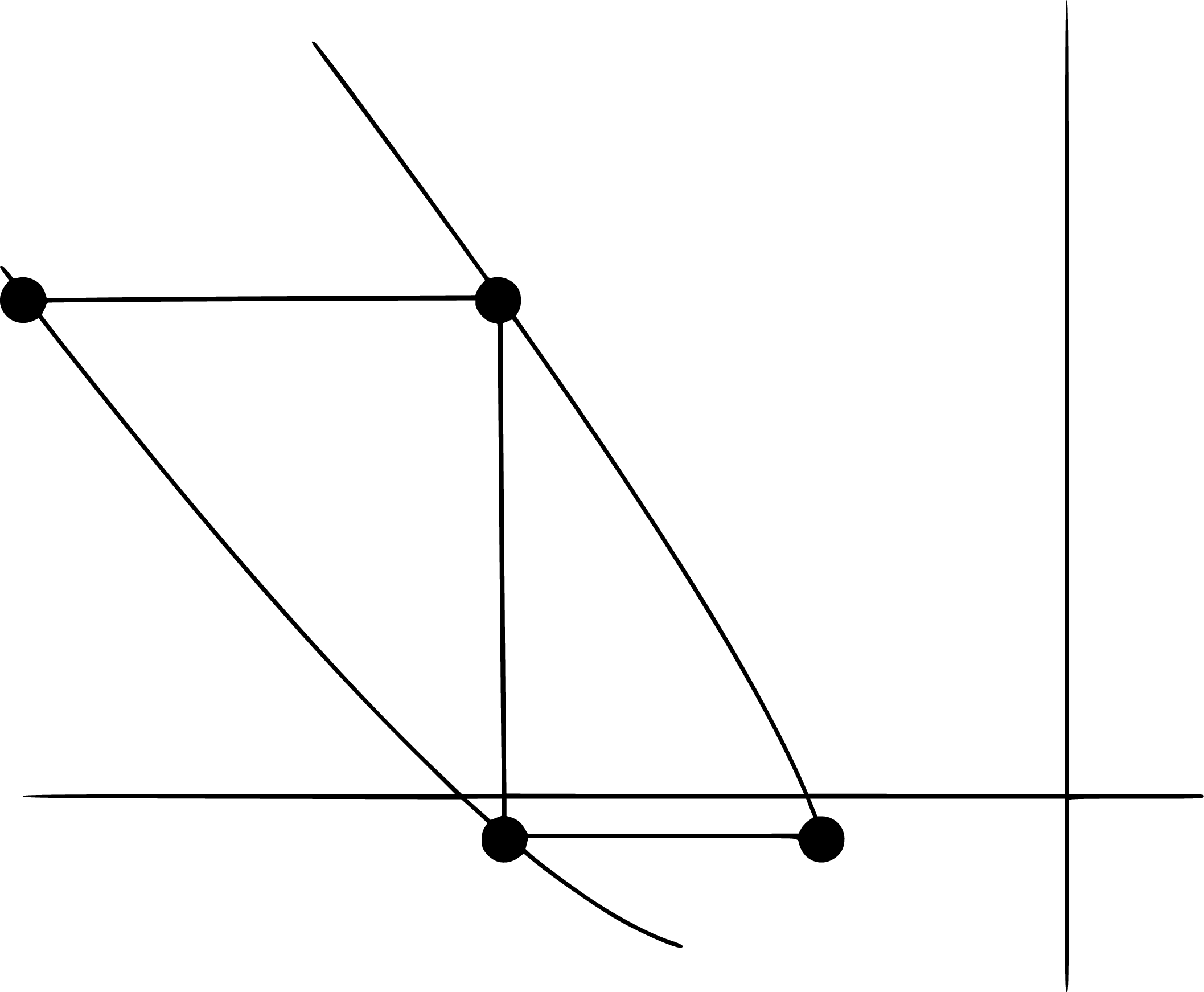}\caption{\small When $|x_2|>M$ or $|x_3|>M$\normalsize}\label{b}
\end{figure}
\end{itemize}
\noindent Therefore, $\mathfrak{U}(\mathbf{k})\subset\{\mathbf{x}\in V_\mathbf{k}(\mathbb{Z}):\exists(s,t,r)\in\mathcal{I}\text{ s.t. }|x_s|=1\ \&\ |x_t|,|x_r|\leq M\}$ is finite.
\end{proof}
\begin{remark}\label{rmk}
    Combining the results of Lemmas \ref{2.1}, \ref{2.2}, and \ref{2.3}, we conclude that each graph $\mathcal{G}$ associated with $\Delta$ has only finitely many minimal vertices. Moreover, if $\mathfrak{U}(\mathbf{k})=\emptyset$, our choice of exactly one representative from each connected cluster of adjacent minimal vertices within each of the sets $\mathfrak{S}(\mathbf{k})$ and $\mathfrak{T}_s(\mathbf{k})$ for $s=0,\pm1$ implies that each $\mathcal{G}$ has exactly one minimal vertex.
\end{remark}


\section{Proofs of Theorems \ref{1.1}, \ref{1.2}, and \ref{1.3}} 
\subsection{\emph{Proof of Theorem \ref{1.1}}} Since every element in $V_\mathbf{k}(\mathbb{Z})$ is $\Gamma_\mathbf{k}$-equivalent to a minimal vertex, it suffices to show that every minimal vertex is $\Gamma_\mathbf{k}$-equivalent to an element in $\mathfrak{S}(\mathbf{k})\cup\mathfrak{T}_0(\mathbf{k})\cup\mathfrak{T}_1(\mathbf{k})\cup\mathfrak{T}_{-1}(\mathbf{k})$.
Note that 
\begin{align*}
    \mathfrak{S}(\mathbf{k})\cup\mathfrak{T}_0(\mathbf{k})\cup\mathfrak{T}_1(\mathbf{k})\cup\mathfrak{T}_{-1}(\mathbf{k})&\subset\{\mathbf{x}\in V_\mathbf{k}(\mathbb{Z}):\forall(s,t,r)\in\mathcal{I},\ |x_s|\leq|\alpha_s-x_tx_r-x_s|\}\\
    &=\{\mathbf{x}\in V_\mathbf{k}(\mathbb{Z}):\text{$\mathbf{x}$ is a minimal vertex}\}.
\end{align*}
Consider $\{\mathbf{x}\in V_\mathbf{k}(\mathbb{Z}):\text{$\mathbf{x}$ is a minimal vertex}\}=\mathfrak{M}_{<2}(\mathbf{k})\cup\mathfrak{M}_{\geq2}(\mathbf{k})$, where
\begin{align*}
    \mathfrak{M}_{<2}(\mathbf{k})&=\{\mathbf{x}\in V_\mathbf{k}(\mathbb{Z}):\text{$\mathbf{x}$ is a minimal vertex}\}\cap\{\mathbf{x}\in V_\mathbf{k}(\mathbb{Z}):\exists s\in\{1,2,3\}\text{ s.t. }|x_s|<2\},\\
    \mathfrak{M}_{\geq2}(\mathbf{k})&=\{\mathbf{x}\in V_\mathbf{k}(\mathbb{Z}):\text{$\mathbf{x}$ is a minimal vertex}\}\cap\{\mathbf{x}\in V_\mathbf{k}(\mathbb{Z}):\forall s\in\{1,2,3\},\ |x_s|\geq2\}.
\end{align*}
Then $\mathfrak{T}_0(\mathbf{k})\cup\mathfrak{T}_1(\mathbf{k})\cup\mathfrak{T}_{-1}(\mathbf{k})\subset\mathfrak{M}_{<2}(\mathbf{k})$ and $\mathfrak{S}(\mathbf{k})\subset\mathfrak{M}_{\geq2}(\mathbf{k})$. By Lemma \ref{2.1}, every element of $\mathfrak{M}_{<2}(\mathbf{k})$ is $\Gamma_\mathbf{k}$-equivalent to one in $\mathfrak{T}_0(\mathbf{k})\cup\mathfrak{T}_1(\mathbf{k})\cup\mathfrak{T}_{-1}(\mathbf{k})$. Suppose $\mathbf{x}\in\mathfrak{M}_{\geq2}(\mathbf{k})$. For all $(s,t,r)\in\mathcal{I}$ with $x_tx_r\neq\alpha_s$, we have $\alpha_s-x_tx_r-x_s\neq-x_s$, that is, $|x_s|<|\alpha_s-x_tx_r-x_s|$. For all $(s,t,r)\in\mathcal{I}$ with $x_tx_r=\alpha_s$, we have $\alpha_s-x_tx_r-x_s=-x_s$, which implies $(\pm x_s,x_t,x_r)$ are $\Gamma_\mathbf{k}$-equivalent, so we may assume $|x_s|=x_s$. In either case, it follows that $\mathbf{x}\in\mathfrak{S}(\mathbf{k})$. Hence, every minimal vertex is $\Gamma_\mathbf{k}$-equivalent to an element in $\mathfrak{S}(\mathbf{k})\cup\mathfrak{T}_0(\mathbf{k})\cup\mathfrak{T}_1(\mathbf{k})\cup\mathfrak{T}_{-1}(\mathbf{k})$.

On the other hand, $\mathfrak{S}(\mathbf{k})\cup\mathfrak{T}_0(\mathbf{k})\cup\mathfrak{T}_1(\mathbf{k})\cup\mathfrak{T}_{-1}(\mathbf{k})$ is finite by Lemma \ref{2.3}. As mentioned in Remark \ref{rmk}, if $\mathfrak{U}(\mathbf{k})=\emptyset$ then each graph associated with $\Delta$ has exactly one minimal vertex, which means all elements in $\mathfrak{S}(\mathbf{k})\cup\mathfrak{T}_0(\mathbf{k})\cup\mathfrak{T}_1(\mathbf{k})\cup\mathfrak{T}_{-1}(\mathbf{k})$ are not $\Gamma_\mathbf{k}$-equivalent to each other.\qed

\subsection{\emph{Proof of Theorem \ref{1.2}}} Except in the case in which a minimal vertex satisfies $-x_i=\alpha_i-x_jx_k-x_i$ (i.e., when two adjacent vertices are both minimal), every $\mathbf{x}\in V_{\mathbf{k}}(\mathbb{Z})$ has the graph corresponding to $\Gamma_{\mathbf{k}}\cdot\mathbf{x}$ that contains exactly one minimal vertex lying in one of the sets $\mathfrak{S}(\mathbf{k})$, $\mathfrak{T}_0(\mathbf{k})$, $\mathfrak{T}_1(\mathbf{k})$, or $\mathfrak{T}_{-1}(\mathbf{k})$, unless the graph contains a vertex in $\mathfrak{U}(\mathbf{k})$. As in Lemma \ref{2.3}, there is $M>0$ such that $|x_1|,|x_2|,|x_3|\leq M$ for all $\mathbf{x}\in\mathfrak{U}(\mathbf{k})$. Thus, all elements of $\mathfrak{U}(\mathbf{k})$ can be found by iterating the Vieta involutions for only finitely many steps. Write $\mathfrak{U}(\mathbf{k})=\{\mathbf{u}_1,\dots,\mathbf{u}_N\}$. For each $l\in\{1,\dots,N\}$, we then collect all minimal vertices whose ascending vertex is $\mathbf{u}_l$ and define
\begin{align*}
    \mathfrak{M}(\mathbf{u}_l)&=\{\mathbf{x}\in V_\mathbf{k}(\mathbb{Z}):\text{$\mathbf{x}$ is a minimal vertex}\}\cap\{\mathbf{x}\in V_\mathbf{k}(\mathbb{Z}):\text{$\mathbf{x}$ is $\Gamma_\mathbf{k}$-equivalent to $\mathbf{u}_l$}\}.
\end{align*}
We then define $\mathfrak{M}_l$, the collection of all minimal vertices that that arise as descending vertices from $\mathfrak{U}(\mathbf{k})$ and are $\Gamma_{\mathbf{k}}$-equivalent to elements of $\mathfrak{M}(\mathbf{u}_l)$, by the following inductive procedure: Set $\mathfrak{M}_l^{(1)}=\mathfrak{M}(\mathbf{u}_l)$ and $s_1=l$. Given $\mathfrak{M}_l^{(t)}$, if $\exists s_{t+1}\in\{1,\dots,N\}\setminus\{s_1,\dots,s_t\}$ such that $\mathfrak{M}_l^{(t)}\cap \mathfrak{M}(\mathbf{u}_{s_{t+1}})\neq\emptyset$, define $\mathfrak{M}_l^{(t+1)}:=\mathfrak{M}_l^{(t)}\cup\mathfrak{M}(\mathbf{u_{s_{t+1}}})$. Otherwise, if $\mathfrak{M}_l^{(t)}\cap \mathfrak{M}(\mathbf{u}_{s})=\emptyset$ $\forall s\in\{1,\dots,N\}\setminus\{s_1,\dots,s_t\}$, then we terminate the process and set $\mathfrak{M}_l=\mathfrak{M}_l^{(t)}$. By construction, this procedure terminates after finitely many steps since $\mathfrak{U}(\mathbf{k})$ is finite and $\Delta(\mathbf{x})$ is $\mathbb{Z}_{\geq0}$-value. Moreover, $\mathfrak{M}(\mathbf{u}_l)=\mathfrak{M}(\mathbf{u}_r)$ if and only if $\mathfrak{M}(\mathbf{u}_l)\cap \mathfrak{M}(\mathbf{u}_r)\neq\emptyset$ (i.e.,\ $\mathbf{u}_l$ and $\mathbf{u}_r$ are $\Gamma_{\mathbf{k}}$-equivalent). 

Given two integral solutions $\mathbf{x},\mathbf{y}\in V_\mathbf{k}(\mathbb{Z})$, let $\mathbf{m}_\mathbf{x}$ and $\mathbf{m}_\mathbf{y}$ denote the minimal vertices whose ascending vertices are $\mathbf{x}$ and $\mathbf{y}$, respectively.
If $\mathbf{x}$ and $\mathbf{y}$ are  $\Gamma_\mathbf{k}$-equivalent, then exactly one of the following occurs: (i) $\mathbf{m}_\mathbf{x}=\mathbf{m}_\mathbf{y}$; (ii) $\mathbf{m}_\mathbf{x},\mathbf{m}_\mathbf{y}$ are adjacent minimal vertices; (iii) $\mathbf{m}_\mathbf{x},\mathbf{m}_\mathbf{y}\in\mathfrak{M}_l$ for some $l\in\{1,\dots,N\}$. Otherwise, $\mathbf{x}$ and $\mathbf{y}$ are not $\Gamma_\mathbf{k}$-equivalent.

Since the computation of $\{\mathbf{u}_1,\dots,\mathbf{u}_N\}$,  $\mathbf{m}_\mathbf{x}$, and $\mathbf{m}_\mathbf{y}$ requires only finitely many iterations of the Vieta involutions, the entire procedure terminates after finitely many steps. \qed


\subsection{\emph{Proof of Theorem \ref{1.3}}} We introduce graphs associated with $\Delta$ corresponding to the mapping class group orbits, including the notions of ascending and descending vertices and minimal vertices, by replacing each Vieta involution $\mathcal{V}_i$ with the composition $\mathcal{V}_i\mathcal{V}_j$. That is, the definitions used for $\Gamma_\mathbf{k}$-orbits are analogously adapted for the mapping class group orbits. 

The mapping class group $\MCG(\Sigma_{0,4})$ is an index-$2$ subgroup of $\Gamma_\mathbf{k}$, generated by compositions of two Vieta involutions. Consequently, the graph corresponding to each $\MCG(\Sigma_{0,4})$-orbit inherits the edge configuration from the corresponding $\Gamma_{\mathbf{k}}$-orbit, where a vertex $\mathbf{z}$ is connected with $\mathcal{V}_2\mathcal{V}_1(\mathbf{z})$, $\mathcal{V}_3\mathcal{V}_2(\mathbf{z})$, and $\mathcal{V}_1\mathcal{V}_3(\mathbf{z})$ by edges. 

If $\mathbf{x}=(x_1,x_2,x_3)\in V_\mathbf{k}(\mathbb{Z})$ with $|x_1|\leq1$, then we have the following (see Section 2.2):
\begin{align*}
    \begin{tikzpicture}[
node distance = 8mm and 8mm,
     V/.style = {},
every edge quotes/.style = {auto, font=\footnotesize}
                    ]
    \begin{scope}[nodes=V]
\node (1)   {\footnotesize{$(0,x_2,x_3)$}};
\node (2) [below = of 1]          {\footnotesize{$(0,\alpha_2-x_2,\alpha_3-x_3)$}};
    \end{scope}
\draw   (2)  edge["\footnotesize{$\mathcal{V}_3\mathcal{V}_2$}"] (1)
        (1)  edge["\footnotesize{$\mathcal{V}_2\mathcal{V}_3$}"] (2);
    \end{tikzpicture}
\hspace{15pt}
    \begin{tikzpicture}[
node distance = 8mm and 8mm,
     V/.style = {},
every edge quotes/.style = {auto, font=\footnotesize}
                    ]
    \begin{scope}[nodes=V]
\node (1)   {\footnotesize{$(\pm1,x_2,x_3)$}};
\node (2) [below left= of 1]    {\footnotesize{$(\pm1,\alpha_2\mp x_3-x_2,\alpha_3\mp(\alpha_2-x_2))$}};
\node (3) [right= of 2]    {\footnotesize{$(\pm1,\alpha_2\mp(\alpha_3-x_3),\alpha_3\mp x_2-x_3)$}};
    \end{scope}
\draw   (2)  edge["\footnotesize{$\mathcal{V}_3\mathcal{V}_2$}"] (1)
        (1)  edge["\footnotesize{$\mathcal{V}_2\mathcal{V}_3$}"] (3)
        (3)  edge["\footnotesize{$\mathcal{V}_2\mathcal{V}_3$}"] (2)
        (2)  edge["\footnotesize{$\mathcal{V}_3\mathcal{V}_2$}"] (3);
    \end{tikzpicture}
\end{align*}
Now, consider $\mathfrak{V}(k)=\{\mathcal{V}_i(\mathbf{x}),\mathcal{V}_j\mathcal{V}_i(\mathbf{x}):\mathbf{x}\in\mathfrak{U}(k)\ \text{and}\ i,j\in\{1,2,3\}\}$. The following two claims are implied by Lemmas \ref{2.3} and \ref{2.2}, respectively.
\begin{claim}\label{3.1}
    $\mathfrak{V}(k)$ is finite.
\end{claim}
\begin{claim}\label{3.2} Let $\mathbf{x}=(x_1,x_2,x_3)\in V_\mathbf{k}(\mathbb{Z})$, and let $\mathcal{G}_\mathbf{x}$ the graph corresponding to its mapping class group orbit. If $|x_s|\geq3$ and $|x_tx_r|>3|\alpha_s|$ for all $(s,t,r)\in\mathcal{I}$ and if $\mathbf{x}$ is a minimal vertex of $\mathcal{G}_\mathbf{x}$, then $\mathcal{G}_\mathbf{x}$ has at most two other minimal vertices; in particular, it has at most three minimal vertices.
\end{claim}
Let us elaborate on the proof of Claim \ref{3.2}: Suppose $\mathbf{x}$ is a minimal vertex of $\mathcal{G}_\mathbf{x}$ with $|x_s|\geq3$, $|x_tx_r|>3|\alpha_s|$ for all $(s,t,r)\in\mathcal{I}$. In the mapping class group orbit, each ascending vertex from $\mathbf{x}$ has exactly one descending vertex as in the proof of Lemma \ref{2.2}. This implies that $\mathcal{G}_\mathbf{x}$ has at most three minimal vertices, which are connected by an edge.
\begin{align*}
    \resizebox{\hsize}{!}{\begin{tikzpicture}[
node distance = 5mm and 8mm,
     V/.style = {},
every edge quotes/.style = {auto, font=\footnotesize}
                    ]
    \begin{scope}[nodes=V]
\node (1)   {\footnotesize{$\mathbf{x}$}};
\node (2) [above left = of 1]    {\footnotesize{$\mathcal{V}_1(\mathbf{x})$}};
\node (3) [above = of 1]          {\footnotesize{$\mathcal{V}_2(\mathbf{x})$}};
\node (4) [above right = of 1]    {\footnotesize{$\mathcal{V}_3(\mathbf{x})$}};
\node (5) [above left = of 2]    {\footnotesize{$\mathcal{V}_3\mathcal{V}_1(\mathbf{x})$}};
\node (6) [left= of 5]          {\footnotesize{$\mathcal{V}_2\mathcal{V}_1(\mathbf{x})$}};
\node (7) [right = of 5]    {\footnotesize{$\mathcal{V}_3\mathcal{V}_2(\mathbf{x})$}};
\node (8) [above = of 3]          {\footnotesize{$\mathcal{V}_1\mathcal{V}_2(\mathbf{x})$}};
\node (9) [right = of 8]    {\footnotesize{$\mathcal{V}_1\mathcal{V}_3(\mathbf{x})$}};
\node (10) [above right= of 4]          {\footnotesize{$\mathcal{V}_2\mathcal{V}_3(\mathbf{x})$}};
\node (11) [above = of 5]          {\footnotesize{$\cdots$}};
\node (12) [above = of 6]          {\footnotesize{$\cdots$}};
\node (13) [above = of 7]          {\footnotesize{$\cdots$}};
\node (14) [above = of 8]          {\footnotesize{$\cdots$}};
\node (15) [above = of 9]          {\footnotesize{$\cdots$}};
\node (16) [above = of 10]          {\footnotesize{$\cdots$}};
    \end{scope}
\draw   (1)  edge["\footnotesize{$\mathcal{V}_1$}"] (2)
        (3)  edge["\footnotesize{$\mathcal{V}_2$}"] (1)
        (4)  edge["\footnotesize{$\mathcal{V}_3$}"] (1)
        (5)  edge["\footnotesize{$\mathcal{V}_3$}"] (2)
        (2)  edge["\footnotesize{$\mathcal{V}_2$}"] (6)
        (3)  edge["\footnotesize{$\mathcal{V}_3$}"] (7)
        (8)  edge["\footnotesize{$\mathcal{V}_1$}"] (3)
        (4)  edge["\footnotesize{$\mathcal{V}_1$}"] (9)
        (10)  edge["\footnotesize{$\mathcal{V}_2$}"] (4)
        (11)  edge[""] (5)
        (12)  edge[""] (6)
        (13)  edge[""] (7)
        (14)  edge[""] (8)
        (15)  edge[""] (9)
        (16)  edge[""] (10);
    \end{tikzpicture}\begin{tikzpicture}[
node distance = 5mm and 8mm,
     V/.style = {},
every edge quotes/.style = {auto, font=\footnotesize}
                    ]
    \begin{scope}[nodes=V]
\node (1)   {\footnotesize{$\mathcal{V}_1(\mathbf{x})$}};
\node (2) [above left = of 1]    {\footnotesize{$\mathbf{x}$}};
\node (3) [above = of 1]          {\footnotesize{$\mathcal{V}_2\mathcal{V}_1(\mathbf{x})$}};
\node (4) [above right = of 1]    {\footnotesize{$\mathcal{V}_3\mathcal{V}_1(\mathbf{x})$}};
\node (5) [above = of 2]    {\footnotesize{$\mathcal{V}_3(\mathbf{x})$}};
\node (6) [above left= of 2]          {\footnotesize{$\mathcal{V}_2(\mathbf{x})$}};
\node (8) [above right = of 3]          {\footnotesize{$\cdots$}};
\node (9) [above right = of 4]    {\footnotesize{$\cdots$}};
\node (11) [above = of 5]          {\footnotesize{$\mathcal{V}_1\mathcal{V}_3(\mathbf{x})$}};
\node (13) [above right = of 5]          {\footnotesize{$\mathcal{V}_2\mathcal{V}_3(\mathbf{x})$}};
\node (12) [above = of 6]          {\footnotesize{$\mathcal{V}_3\mathcal{V}_2(\mathbf{x})$}};
\node (14) [above left = of 6]          {\footnotesize{$\mathcal{V}_1\mathcal{V}_2(\mathbf{x})$}};
\node (16) [above = of 11]          {\footnotesize{$\cdots$}};
\node (17) [above = of 12]          {\footnotesize{$\cdots$}};
\node (18) [above = of 13]          {\footnotesize{$\cdots$}};
\node (19) [above = of 14]          {\footnotesize{$\cdots$}};
    \end{scope}
\draw   (1)  edge["\footnotesize{$\mathcal{V}_1$}"] (2)
        (3)  edge["\footnotesize{$\mathcal{V}_2$}"] (1)
        (4)  edge["\footnotesize{$\mathcal{V}_3$}"] (1)
        (5)  edge["\footnotesize{$\mathcal{V}_3$}"] (2)
        (2)  edge["\footnotesize{$\mathcal{V}_2$}"] (6)
        (8)  edge[""] (3)
        (4)  edge[""] (9)
        (11)  edge["\footnotesize{$\mathcal{V}_1$}"] (5)
        (13)  edge["\footnotesize{$\mathcal{V}_2$}"] (5)
        (12)  edge["\footnotesize{$\mathcal{V}_3$}"] (6)
        (6)  edge["\footnotesize{$\mathcal{V}_1$}"] (14)
        (16)  edge[""] (11)
        (17)  edge[""] (12)
        (18)  edge[""] (13)
        (19)  edge[""] (14);
    \end{tikzpicture}}
    \end{align*}
    
Combining Claims \ref{3.1} and \ref{3.2}, we conclude that each mapping class group orbit has at most finitely many minimal vertices.

\vspace{5pt}

\emph{Proof of Theorem \ref{1.3}.} Let $\mathbf{x},\mathbf{y}\in V_k(\mathbb{Z})$ be given. Since $\mathfrak{V}(\mathbf{k})$ is finite by Claim \ref{3.1}, we can apply the same method as in the proof of Theorem \ref{1.2}.\qed 

\begin{remark}
    In all proofs in this paper, $\beta$ plays no role. Therefore, the results remain valid even without imposing any relation between the $\alpha_i$ and $\beta$. In the special case $(x,y,z)=(-x_1,x_2,x_3)$, with $\alpha_1=\alpha_2=\alpha_3=0$ and $\beta=k$, the results specialize to those for the character variety of a one-holed torus (see example \ref{ex:1.5}).
\end{remark}

\bibliographystyle{plain}
\bibliography{bib}

\end{document}